\documentclass[10pt]{amsart}
\usepackage{amssymb}

\newtheorem{theorem}{Theorem}
\newtheorem{lemma}[theorem]{Lemma}

\newtheorem{corollary}[theorem]{Corollary}

\newcommand{\newsection}[1] {\section{#1}\setcounter{theorem}{0}
 \setcounter{equation}{0}\par\noindent}

\newcommand{\cal}{\mathcal}

\newcommand{\supp}{\text{supp}}
\newcommand{\dist}{\text{dist}}
\newcommand{\eps}{\varepsilon}
\newcommand{\Rn}{{{\mathbb R}^n}}
\newcommand{\R}{{\mathbb R}}

\newcommand{\g}{{\rm g}}
\newcommand{\cd}{\,\cdot\,}

\newcommand{\Dy}{D_{\!y}}

\renewcommand{\l}{\lambda}
\newcommand{\demo}{\noindent{\it Proof.} }
\newcommand{\findemo}{\qed\bigskip}
\newcommand{\thetabar}{{\bar{\theta}}}
\newcommand{\U}{{\cal U}}

\begin{document}

\title[$L^p$ norm of Spectral clusters]
{On the $L^p$ norm of spectral clusters for compact manifolds with
boundary}
\thanks{The authors were supported by the National Science Foundation,
Grants DMS-0140499, DMS-0099642, and DMS-0354668.}
\author{Hart F. Smith}
\author{Christopher D. Sogge}
\address{Department of Mathematics, University of Washington,
Seattle, WA 98195}
\address{Department of Mathematics, Johns Hopkins University,
Baltimore, MD 21218}
\email{hart@math.washington.edu}
\email{sogge@jhu.edu}

\maketitle

\newsection{Introduction}\label{section1}
Let $M$ be a compact two-dimensional manifold with boundary, and
let $P$ be an elliptic, second order differential operator on $M$,
self-adjoint with respect to a density $d\mu$, and with
vanishing zeroeth order term, so that in local coordinates
\begin{equation}\label{Pform}
\bigl(Pf\bigr)(x)=
\rho(x)^{-1}\sum_{i,j=1}^n
\partial_i\Bigl(\rho(x)\,\g^{ij}(x)\,\partial_jf(x)\Bigr)
\,,\qquad d\mu=\rho(x)\,dx\,.
\end{equation}
We take $\g^{ij}$ to be positive, so that the Dirichlet eigenvalues of
$P$ can be written as $\{-\lambda_j^2\}_{j=0}^\infty\,.$

Let $\chi_\l$ be the projection of $L^2(d\mu)$ onto the subspace
spanned by the Dirichlet eigenfunctions for which
$\lambda_j\in[\l,\l+1]$. In the case that $M$ is compact without
boundary of dimension $n\ge 2$, and the coefficients of $P$ are
$C^\infty$ functions, Sogge \cite{So} established the following
bounds
\begin{equation}\label{estimate1a}
\bigl\|\chi_\l f\bigr\|_{L^q(M)}\le C\,\l^{\frac{n-1}2(\frac
12-\frac 1q)}\,\|f\|_{L^2(M)}\,, \qquad 2\le q\le q_n\,.
\end{equation}
\begin{equation}\label{estimate2a}
\bigl\|\chi_\l f\bigr\|_{L^q(M)}\le C\,\l^{n(\frac 12-\frac
1q)-\frac 12}\,\|f\|_{L^2(M)}\,, \qquad q_n\le q\le\infty\,.
\end{equation}
Furthermore, the exponent of $\l$ is sharp on every such manifold
(see e.g., \cite{Soggebook}). In the case of a sphere, the
examples which prove sharpness are in fact eigenfunctions. For
\eqref{estimate1a} the appropriate example is an eigenfunction
which concentrates in a $\lambda^{-\frac 12}$ diameter tube about
a geodesic. For \eqref{estimate2a}, the example is a zonal
eigenfunction of $L^2$ norm $\l^{\frac{n-1}2}$ which takes on
value comparable to $\l$ on a $\l^{-1}$ diameter ball about each
of the north and south poles. Approximate spectral clusters with
similar properties can be constructed in the interior of any
smooth manifold, showing that for spectral clusters (though not
necessarily eigenfunctions) the exponents in \eqref{estimate1a}
and \eqref{estimate2a} are also lower bounds on manifolds with
boundary.

In \cite{SmSo}, the authors showed that, on a manifold of dimension $n\ge 2$
for which the boundary is everywhere strictly geodesically concave
(such as the complement in $\Rn$ of a strictly convex set) the estimates
\eqref{estimate1a} and \eqref{estimate2a} both hold.

On the other hand, Grieser \cite{Gr} observed that in the unit
disk $\{|x|\le 1\}$ there are eigenfunctions of the Laplacian, for
Dirichlet as well as for Neumann boundary conditions, of
eigenvalue $-\l^2$ that concentrate within a $\l^{-\frac 23}$
neighborhood of the boundary. These are the classical Rayleigh
whispering gallery modes (see \cite{Ray1}, \cite{Ray2}). The
Fourier-Airy calculus of Melrose and Taylor allows one to
construct an approximate spectral cluster with similar
localization properties near any boundary point of $M$ at which
the boundary is strictly convex (the gliding case). Consequently,
if $M$ is of dimension two and the boundary has a point of strict
convexity with respect to the metric $\g$ (for instance, any
smoothly bounded planar domain endowed with the standard Laplacian and
either Dirichlet or Neumann conditions) the following bounds
cannot be improved upon
\begin{equation}\label{estimate1}
\bigl\|\chi_\l f\bigr\|_{L^q(M)}\le C\,\l^{\frac 23(\frac 12-\frac
1q)}\,\|f\|_{L^2(M)}\,, \qquad 2\le q\le 8\,.
\end{equation}
\begin{equation}\label{estimate2}
\bigl\|\chi_\l f\bigr\|_{L^q(M)}\le C\,\l^{2(\frac 12-\frac
1q)-\frac 12}\,\|f\|_{L^2(M)}\,, \qquad 8\le q\le\infty\,.
\end{equation}

In this paper we show that the estimates \eqref{estimate1}
and \eqref{estimate2} hold on any two dimensional compact manifold
with boundary, for $P$ as above and either Dirichlet or Neumann
conditions assumed.
Estimate \eqref{estimate1} follows
by interpolation of the trivial case $q=2$ with the case $q=6$,
so we restrict attention to $q\ge 6$ for \eqref{estimate1}.
For $q\ge 6$, the above estimates are an immediate consequence
of the following theorem (see for example \cite{MSS} or \cite{Sm1}).

\begin{theorem}\label{maintheorem1}
Suppose that $u$ solves the Cauchy problem on $\R\times M$
\begin{equation}\label{cauchy}
\partial_t^2u(t,x)=Pu(t,x)\,,
\qquad
u(0,x)=f(x)\, ,
\qquad \partial_tu(0,x)=0\,,
\end{equation}
and satisfies either Dirichlet conditions
\begin{equation*}
 u(t,x)=0\quad\text{if}\quad x\in\partial M\,,
\end{equation*}
or Neumann conditions, where $N_x$ is a
unit normal field with respect to $\g$,
\begin{equation*}
N_x\cdot\nabla_x u(t,x)=0\quad\text{if}\quad x\in\partial M\,.
\end{equation*}
Then the following bounds hold for $6\le q \le 8$,
\begin{equation*}
\|u\|_{L^q_xL^2_t(M\times[-1,1])}
\le C\,\|f\|_{H^{\gamma(q)}(M)}\,,
\qquad \textstyle
\gamma(q)=\frac23\bigl(\frac 12-\frac 1q\bigr)\,,
\end{equation*}
and the following bounds hold for $8\le q\le \infty$,
\begin{equation*}
\|u\|_{L^q_xL^2_t(M\times[-1,1])}
\le C\,\|f\|_{H^{\delta(q)}(M)}\,.
\qquad \textstyle
\delta(q)=2\bigl(\frac 12-\frac 1q\bigr)-\frac 12\,.
\end{equation*}
\end{theorem}

In the statement of the theorem, the space $H^s(M)$ refers to the
Sobolev space of order $s$ on $M$ determined, respectively,
by Dirichlet or Neumann eigenfunctions.

Our approach to proving Theorem \ref{maintheorem1} is to work
in geodesic normal coordinates near $\partial M$, and to extend
both the operator $P$ and the solution $u$ across the boundary,
to obtain $u$ as a solution to a wave equation on an open set,
but for an operator with coefficients of Lipschitz regularity.
We then adapt a frequency dependent scaling argument, originally developed
to handle Lipschitz metrics, to metrics with the particular type
of codimension-1 singularities that the extended $P$ will have.

We remark that, for operators of the type \eqref{Pform} with
$\rho$ and $\g^{ij}$ of Lipschitz regularity, the estimate
\eqref{estimate1} is known on the range $2\le q\le 6$, as
established by the first author in \cite{Sm2}, along with a weaker
version of \eqref{estimate2} having larger exponent if $q<\infty$.
It is not currently known what the sharp exponents are for general
Lipschitz $P$, since the known counterexamples satisfy the
estimates \eqref{estimate2}.  The estimates for $q=\infty$ were
established for eigenfunctions recently by Grieser \cite{Gr2},
while the sup-norm estimates for spectral clusters were obtained
by the second author in \cite{So2}.

For $q=\infty$, the squarefunction estimate of Theorem \ref{maintheorem2}
below was shown in \cite{Sm2} to hold for operators $P$ with Lipschitz
coefficients, which in particular
implies the $q=\infty$ case of Theorem \ref{maintheorem1}
for $P$ on a manifold with boundary. Our proof here
of the case $q<\infty$, however, depends crucially on the fact that
if $u$ is appropriately microlocalized away
from directions tangent to $\partial M$,
then better squarefunction estimates hold
than do for directions near to tangent.
In other words, we exploit the fact that
the more highly localized eigenfunctions
considered in \cite{Gr} are associated only
to gliding directions along $\partial M$,
not directions transverse to $\partial M$.

A historical curiosity is that the critical $L^2\to L^8$ bounds
for $\chi_\lambda$ have an analog in Euclidean space which seems
to be the first restriction theorem for the Fourier transform.  To
explain this, we first notice that by duality our $L^2\to L^8$
bounds are equivalent to the statement that $\chi_\lambda:
L^{8/7}\to L^2$ with norm $O(\lambda^{1/4})$.  The Euclidean analog
would say that if $\chi_\lambda: L^{8/7}({\mathbb R}^2)\to
L^2({\mathbb R}^2)$ denotes the projection onto Fourier
frequencies $|\xi|\in [\lambda,\lambda+1]$, then this operator
also has norm $O(\lambda^{1/4})$.  An easy scaling argument shows
then that the latter result is equivalent to the following Fourier
restriction theorem for the circle
$$\Bigl(\int_0^{2\pi}|\Hat f(\cos\theta,\sin\theta)|^2\,
d\theta \Bigr)^{1/2}\le
C\|f\|_{L^{8/7}({\mathbb R}^2)}, \quad f\in C^\infty_0({\mathbb
R}^2).$$ Stein \cite{St} proved this using a now standard $TT^*$
argument, together with DeLeeuw's \cite{DeL} observation that
$\widehat{d\theta}$ maps $L^{8/7}({\mathbb R}^2)\to L^8({\mathbb
R}^2)$ by the Hardy-Littlewood-Sobolev theorem, as
$|\widehat{d\theta}|\le C|x|^{-1/2}$.  Since this argument does
not use the oscillations of $\widehat{d\theta}$, one can
strengthen the above restriction theorem to show that, for $j\ge 1$,
one has the uniform bounds
\begin{equation}\label{1.7}
\Bigl(\, \int_0^{2^{-j}}
|\hat f(\cos\theta,\sin\theta)|^2\, d\theta\,
\Bigr)^{1/2}\le C\,2^{-j/8}\|f\|_{L^{8/7}({\mathbb R}^2)}, \quad
f\in C^\infty_0({\mathbb R}^2).
\end{equation}
By the Knapp example, there is no small angle improvement
for the critical $L^{6/5}({\mathbb R}^2)\to L^2({\mathbb
S}^1)$ restriction theorem of Stein-Tomas.  A key step for us is
that in the setting of compact manifolds with boundary we also get
the same $O(2^{-j/8})$ improvement in our $L^8$-estimates when
microlocalized to regions of phase space that correspond to
bicharacteristics that are of angle comparable to $2^{-j}$ from
tangency to the boundary.

In higher dimensions the natural analog of
\eqref{estimate1}-\eqref{estimate2} would say that
\begin{equation}\label{1.8}
\bigl\|\chi_\l f\bigr\|_{L^q(M)}\le C\,\l^{(\frac 23
+\frac{n-2}2)(\frac 12-\frac 1q)}\,\|f\|_{L^2(M)}\,, \qquad 2\le
q\le \tfrac{6n+4}{3n-4}\,.
\end{equation}
\begin{equation}\label{1.9}
\bigl\|\chi_\l f\bigr\|_{L^q(M)}\le C\,\l^{n(\frac 12-\frac
1q)-\frac 12}\,\|f\|_{L^2(M)}\,, \qquad \tfrac{6n+4}{3n-4}\le
q\le\infty\,.
\end{equation}
By higher dimensional versions of the Rayleigh whispering gallery
modes, this would be sharp if true.  At present we are unable to
prove this estimate but, as we shall indicate in the final
section, we can prove the bounds in \eqref{1.9} for the smaller
range of exponents $q\ge 4$ if $n\ge 4$, and $q\ge 5$ if $n=3$.
We hope to return to the
problem of proving sharp results in higher dimensions in a future
work.

\noindent{\bf Notation.} We use the following notation. The symbol
$a\lesssim b$ means that $a\le C\,b$, where $C$ is a constant that
depends only on globally fixed parameters 
(or on $N,\,\alpha,\,\beta$ in case of
inequalities involving general integers.)

For convenience we will let $x_3$ serve as substitute for the time variable
$t$. We use $d=(d_1,d_2,d_3)$ to denote the gradient operator, and
$D=-id$.

\newsection{Dyadic Localization Arguments}\label{section2}
The estimates of Theorem \ref{maintheorem1} hold if $u$ is supported
away from $\partial M$ by the results of \cite{MSS}. Consequently,
by finite propagation velocity and the
use of a smooth partition of unity we may assume that, for $T$ small,
the solution $u(t,x)$ in Theorem \ref{maintheorem1}
is for $|t|\le T$ supported in a suitably small
coordinate patch centered on the boundary. Note that if we establish
Theorem \ref{maintheorem1} on the set $|t|\le T$ for some small $T$,
it then holds for $T=1$ by energy conservation.

We work in boundary normal coordinates for the Riemannian metric $\g_{ij}$
that is dual to $\g^{ij}$ of \eqref{Pform}.
Thus, $x_2>0$ will define the manifold $M$, and $x_1$ is a coordinate
function  on $\partial M$ which we choose
so that $\partial_{x_1}$ is of unit length along $\partial M$.
In these coordinates,
\begin{equation}\label{cond1}
\g_{22}(x_1,x_2)=1\,,\qquad \g_{11}(x_1,0)=1\,,\qquad
\g_{12}(x_1,x_2)=\g_{21}(x_1,x_2)=0\,.
\end{equation}
The metric $\g^{ij}$ for $P$ is the inverse to $\g_{ij}$, and the
same equalities hold for it.

We now extend the coefficient $\g^{11}$ and $\rho$ in an even manner across
the boundary, so that
\begin{equation}\label{cond2}
\g^{11}(x_1,-x_2)=\g^{11}(x_1,x_2)\,,\qquad
\rho(x_1,-x_2)=\rho(x_1,x_2)\,.
\end{equation}
The extended functions are then piecewise smooth,
and of Lipschitz regularity across $x_2=0$. Because $\g$ is diagonal,
the operator $P$ is preserved under the reflection $x_2\rightarrow -x_2.$

After multiplying $\rho(x)$ by a constant, and rescaling variables
if necessary, we may assume that on the ball $|x|<1$ the function
$\rho(x)$ is $C^N({\mathbb R}^2_+)$ close to the function 1, and
$\g^{ij}(x)$ is $C^N({\mathbb R}^2_+)$ close to the euclidean
metric, where $N$ is suitably large, and $c_0$ will be taken
suitably small,
\begin{equation}\label{cond3}
\|\,\rho-1\|_{C^N(\R^2_+)}\le c_0\,,\qquad
\|\,\g^{ij}-\delta^{ij}\|_{C^N(\R^2_+)}\le c_0\, .
\end{equation}
We may then
extend $\rho$ and $\g^{ij}$ globally, preserving conditions
\eqref{cond1}--\eqref{cond3},
so that $P$ is defined globally on $\R^2$ and such that
\begin{equation}\label{cond4}
\rho(x)=1 \,,\qquad
\g^{ij}(x)=\delta^{ij}\quad\text{for}\;\;|x|\ge \frac 34\,.
\end{equation}

We then extend the initial data $f$ and the solution $u$ to be odd in
$x_2$ (respectively even in $x_2$ in case of Neumann conditions).
This extension map is seen to map the Dirichlet (respectively Neumann)
Sobolev space $H^2(\R^n_+)$ to $H^2(\R^n)$, hence $H^\delta(\R^n_+)$ to
$H^\delta(\R^n)$ for $0\le \delta\le 2$. The extended solution $u$ thus solves
the extended equation $\partial_t^2u=Pu$ on $\R\times\R^2$, with the
extended initial data $f$.
The result of Theorem \ref{maintheorem1} is thus a direct
consequence of the following

\begin{theorem}\label{maintheorem2}
Suppose that the operator $P$ takes the form \eqref{Pform}, and
that $\rho$ and $\g$ satisfy conditions \eqref{cond1}--\eqref{cond4} above.
Let $u$ solve the Cauchy problem on $\R\times\R^2$
\begin{equation}\label{cauchy'}
\partial_t^2u(t,x)=Pu(t,x)\,,
\qquad
u(0,x)=f(x)\, ,
\qquad \partial_tu(0,x)=g(x)\,.
\end{equation}
Then the following bounds hold for $6\le q \le 8$,
\begin{equation*}
\|u\|_{L^q_xL^2_t(\R^2\times[-1,1])}
\lesssim \,\bigl(\,\|f\|_{H^{\gamma(q)}}+\|g\|_{H^{\gamma(q)-1}}\bigr)\,,
\qquad \textstyle
\gamma(q)=\frac23\bigl(\frac 12-\frac 1q\bigr)\,,
\end{equation*}
and the following bounds hold for $8\le q\le \infty$,
\begin{equation*}
\|u\|_{L^q_xL^2_t(\R^2\times[-1,1])}
\lesssim\,\bigl(\,\|f\|_{H^{\delta(q)}}+\|g\|_{H^{\delta(q)-1}}\bigr)\,,
\qquad \textstyle
\delta(q)=2\bigl(\frac 12-\frac 1q\bigr)-\frac 12\,.
\end{equation*}
\end{theorem}

We begin by reducing matters to compactly supported $u$ satisfying
an inhomogeneous equation.
Henceforth, we will use notation $x_3=t$.
Let $\phi(x)$ be a smooth even function on $\R^3$, 
equal to 1 for $|x|\le 3/2$,
and vanishing for $|x|\ge 2$.
We may then write
$$
\sum_{j=1}^3 D_j\Bigl(\,a^{jj}(x)\,D_j
(\phi u)(x)\,\Bigr)=\sum_{j=1}^3\, D_j F_j(x)\,,
$$
where
$$
a^{33}(x)=\rho(x)\,,\qquad
a^{jj}(x)=-\rho(x)\,\g^{jj}(x)\quad\text{for}\;\;j=1,2.
$$
We express this equation concisely as $DAD(\phi u)=DF$,
and observe that for $0\le \delta\le 2$
$$
\|\phi u\|_{H^\delta(\R^3)}+\|F\|_{H^\delta(\R^3)}
\;\lesssim\; \|f\|_{H^\delta}+\|g\|_{H^{\delta-1}}\,.
$$
This is a consequence of energy estimates, which hold
separately on $\R^3_+$ and $\R^3_-$,
together with the fact that $DAD(\phi u)$ is compactly
supported and has integral 0, so may be written as $DF$.

We may thus assume that $u(x)$ is supported in the ball $|x|\le 2$,
and need to show that
\begin{align}\label{uesta}
\|u\|_{L^qL^2}\lesssim &
\;\|u\|_{H^{\gamma(q)}}+\|F\|_{H^{\gamma(q)}}\,,\qquad
6\le q \le 8\,,\\
\notag\\
\label{uestb}
\|u\|_{L^qL^2}\lesssim &
\;\|u\|_{H^{\delta(q)}}+\|F\|_{H^{\delta(q)}}\,,\qquad
8\le q\le \infty\,,
\end{align}
where $DADu=DF\,.$

Next let $\Gamma(\xi)$ be a multiplier of order 0, supported in
the set $\frac 14\,|\xi_3|\le|\xi_1,\xi_2|\le 4|\xi_3|\,,$ which equals
1 on the set $\frac 12\,|\xi_3|\le|\xi_1,\xi_2|\le 2|\xi_3|\,.$
The operator $DAD$ is elliptic on the support of $1-\Gamma$, and
we may write
$$
DAD\bigl(1-\Gamma(D)\bigr)u=
\bigl(1-\Gamma(D)\bigr)DF-D\bigl[A,\Gamma(D)\bigr]Du\,.
$$
As a consequence of the Coifman-Meyer commutator theorem \cite{CM} (see
also Proposition 3.6.B of \cite{Tay})
the operator $[A,\Gamma(D)]$ maps $H^{\delta-1}\rightarrow H^\delta$
for $0\le \delta\le 1$. Hence, the right hand side of the above belongs
to $H^{\delta-1}$, and by Sobolev embedding and elliptic regularity
(see, for example, Theorem 2.2.B of \cite{Tay}, which applies in the Sobolev
setting) we have
$$
\|\bigl(1-\Gamma(D)\bigr)u\|_{L^qL^2}\lesssim\|
\bigl(1-\Gamma(D)\bigr) u\|_{H^{\delta(q)+1}} \lesssim
\|u\|_{H^{\delta(q)}}+\|F\|_{H^{\delta(q)}}\,.
$$
Indeed, there is an extra $\frac 12$ derivative in $\delta(q)+1$
beyond the Sobolev index $n(\frac 12-\frac 1q)$, 
so this holds for all $2\le q\le\infty$.
Since $\gamma(q)\ge \delta(q)$ for $q\le 8$, this implies that
\eqref{uesta} and \eqref{uestb} hold for $u$ replaced by $(1-\Gamma(D))u$
on the left hand side.

It thus remains to establish \eqref{uesta} and \eqref{uestb}
with $u$ replaced on the
left by $\Gamma(D)u$.
We take a Littlewood-Paley decomposition in $\xi$ to write
$$
\Gamma(D)u=\sum_{k=1}^\infty \Gamma_k(D)u=\sum_{k=1}^\infty u_k\,,
$$
with $\widehat{u_k}$ is supported in a region where
$|\xi_1,\xi_2|\approx|\xi_3|$
and $|\xi|\approx 2^k$. Since these regions have finite overlap
in the $\xi_3$ axis, we have
$$
\|\Gamma(D)u\|_{L^qL^2}\lesssim \|u_k\|_{L^q\ell^2_k L^2}\lesssim
\|u_k\|_{\ell^2_k L^q L^2}\,,
$$
where we use $q\ge 2$ at the last step.

Now let $A_k$ denote the matrix of coefficients obtained by truncating
the frequencies of $a^{ii}(x)$ to $|\xi|\le c\,2^k$ for a fixed small
$c$.
We then have $DA_kDu_k=DF_k$, where
\begin{equation}\label{Fkform}
F_k=\Gamma_k(D)F+\bigl[A,\Gamma_k(D)\bigr]Du+\bigl(A_k-A)D u_k\,.
\end{equation}
Note that the inhomogeneity
$F_k$ is now localized in frequency to $|\xi_3|\approx|\xi|\approx2^k\,,$
by the frequency localizations of $A_k$ and $u_k$.

We claim that, for $0\le \delta\le 1$,
$$
\sum_{k=1}^\infty 2^{2k\delta}\|F_k\|^2_{L^2}\lesssim
\|u\|_{H^\delta}^2+\|F\|^2_{H^\delta}
\,.
$$
This follows by orthogonality for the first term on the 
right of \eqref{Fkform},
and the last term is handled by the bound 
$\|A-A_k\|_{L^\infty}\lesssim 2^{-k}$.
The middle term is handled by the Coifman-Meyer commutator theorem,
which yields that $\sum_{k=1}^\infty \varepsilon_k
\bigl[A,\Gamma_k(D)\bigr]$ maps $H^{\delta-1}\rightarrow H^\delta$
for all sequences $\varepsilon_k=\pm 1$.

We thus are reduced to establishing
uniform estimates for each dyadically localized piece $u_k$.
We thus fix a frequency scale $\l=2^k$ for the rest of this paper.
We then need to prove the following estimates, where we now set
$DA_\l Du_\l=F_\l$,
\begin{align*}
\|u_\l\|_{L^qL^2(\R^3)}\lesssim &\;
\l^{\gamma(q)}\bigl(\,\|u_\l\|_{L^2(\R^3)}+\l^{-1}\|F_\l\|_{L^2(\R^3)}
\,\bigr)\,,\quad 6\le q\le 8\,,
\\
\\
\|u_\l\|_{L^qL^2(\R^3)}\lesssim &\;
\l^{\delta(q)}\bigl(\,\|u_\l\|_{L^2(\R^3)}+\l^{-1}\|F_\l\|_{L^2(\R^3)}
\,\bigr)\,,\quad 8\le q\le\infty\,.
\end{align*}
Since we are using $x_3$ orthogonality to make this reduction, we
must control the norms of the $u_\l$ globally. However, since $u$ is
supported in the ball of radius 2, it is easy to see that the
norm of $u_\l$ over $|x|\ge 3$ is bounded by $\l^{-1}\|u\|_{L^2}$,
so in fact it suffices to establish the above estimate with
the left hand side norm taken over the cube of sidelength 3.

If we let $v_\l$ denote the localization of $u_\l$ to frequencies
where $|\xi_2|\ge \frac 18\,|\xi_3|$, then the
square function estimates hold for $v_\l$ as on an open manifold,
$$
\|v_\l\|_{L^qL^2(\R^3)}\lesssim\;
\l^{\delta(q)}\bigl(\,\|u_\l\|_{L^2(\R^3)}+\|F_\l\|_{L^2(\R^3)}
\,\bigr)\,,\quad 6\le q\le\infty\,.
$$
This will follow as a consequence of the techniques we use to handle
the part of $u_\l$ with frequencies localized to angle $\approx 1$
from the $\xi_3$ axis.

Consequently, we will assume that
$$
\text{supp}\bigl(\widehat{u_\l}\bigr)\,\subseteq\,
\bigl\{\,\xi\,:\,
|\xi_1|\in\bigl[\tfrac 12\l,2\l\bigr]\,,\quad
|\xi_2|\le\tfrac 1{10}\l\,,
\quad|\xi_3|\in\bigl[\tfrac 12\l,2\l\bigr]\,
\bigr\}\,.
$$
On this region, the operator $DA_\l D$ is hyperbolic with respect
to the $x_1$ direction. We can thus take $p(x,\xi')$ a
positive elliptic symbol
in $\xi'=(\xi_2,\xi_3)$, so that
$$
a^{11}_\l(x)\bigl(\,\xi_1^2-p(x,\xi')^2\,\bigr)
=\sum_{j=1}^3 a^{jj}_\l(x)\xi_j^2
\quad\text{if}\quad
|\xi_2|\le\tfrac 19\l\,,\quad|\xi_3|\in\bigl[\tfrac 13\l,3\l\bigr]\,,
$$
and such that
$$
p(x,\xi')=|\xi'|\quad\text{if}\quad
\xi'\notin\bigl[-\tfrac 18 \l,\tfrac 18\l\bigr]\times
\bigl[\tfrac 14\l,4\l\bigr]\,.
$$
We also smoothly set $p(x,\xi')=1$ near $\xi'=0$.
Thus,
$$
p(x,\xi')\,,\;d_x p(x,\xi')\in S^1_{1,1}\,,
$$
and $p(x,\xi')$ differs from $|\xi'|$ by
a symbol supported in the dyadic shell $|\xi'|\approx\l$.

Next, let $p_\l(x',\xi)$ be obtained by truncating the symbol
$p(x,\xi')$ to $x'$-frequencies less than $c\l$, where $c$ is
a small constant. Then, uniformly over $\l$,
$$
p_\l(x,\xi')-p(x,\xi')\in S^0_{1,1}\,,\quad
\text{support}\,(\,p_\l-p)\subset
\,\bigl\{\,\xi'\,:\,|\,\xi'|\approx\l\,\bigr\}\,.
$$
Furthermore the symbol-composition rule holds for $p_\l$ to first order.
Consequently, we can write
$$
\bigl(\,D_1+p_\l(x,D')\,\bigr)\bigl(\,D_1-p_\l(x,D')\,\bigr)u_\l=F'_\l\,,
$$
where
$$
\|F'_\l\|_{L^2(\R^3)}\lesssim \l\,\|u_\l\|_{L^2(\R^3)}+\|F_\l\|_{L^2(\R^3)}\,.
$$
The function $u_\l$ can be written as the sum of four pieces with disjoint
Fourier transforms, according to the possible signs of $\xi_1$ and $\xi_3$.
We restrict attention to the piece $u_\l^+$,
supported where $\xi_1>0$ and $\xi_3>0$. Estimates for the other pieces
will follow similarly.
Since $p_\l$ is $x'$-frequency localized, $F'_\l$ also splits into
four disjoint pieces. The symbol $\xi_1+p(x,\xi')$ is elliptic on the
region $\xi_1>0$, hence we may write
$$
D_1u_\l^+ - p_\l(x,D')u_\l^+=F''_\l\,,
$$
where
$$
\|F''_\l\|_{L^2(\R^3)}\lesssim
\|u_\l\|_{L^2(\R^3)}+\l^{-1}\|F_\l\|_{L^2(\R^3)}
\,.
$$
Finally, we have that
$$
p_\l(x,D')-p_\l(x,D')^*\in \text{Op}(S^0_{1,1})\,,
$$
and is dyadically supported in $\xi'$.
We have thus reduced the proof of Theorem
\ref{maintheorem2} to the following.
\begin{theorem}\label{maintheorem3}
Suppose that the $x'$-Fourier transform of $u_\l$
satisfies the support condition
$$
{\rm supp}\bigl(\widehat{u_\l}\bigr)\,\subseteq\,
\bigl\{\,\xi'\,:\,
|\xi_2|\le\tfrac 1{10}\l\,,
\quad\xi_3\in\bigl[\tfrac 12\l,2\l\bigr]\,
\bigr\}\,,
$$
and that
$$
D_1 u_\l-P_\l u_\l=F_\l\,,
$$
where $P_\l=\frac 12\bigl(\,p_\l(x,D')+p_\l(x,D')^*\bigr)\,.$
Then, for $S=[0,1]\times\R^2$,
\begin{align}
\|u_\l\|_{L^qL^2(S)}\lesssim &\;
\l^{\gamma(q)}\bigl(\,\|u_\l\|_{L^\infty L^2(S)}+\|F_\l\|_{L^2(S)}
\,\bigr)\,,\quad 6\le q\le 8\,,
\notag\\
\label{uestc}\\
\|u_\l\|_{L^qL^2(S)}\lesssim &\;
\l^{\delta(q)}\bigl(\,\|u_\l\|_{L^\infty L^2(S)}+\|F_\l\|_{L^2(S)}
\,\bigr)\,,\quad 8\le q\le\infty\,.
\notag
\end{align}
\end{theorem}
The use of the $L^\infty L^2$ norm of $u_\l$ is allowed by Duhamel
and energy bounds.  Here, as in what follows, we are using the
shorthand mixed-norm notation that $L^pL^q=L^p_{x_1}L^q_{x'}$.

\newsection{The Angular Localization}\label{section3}
In this section
we take a further decomposition of $u_\l$, by decomposing its
Fourier transform dyadically in the $\xi_2$ variable.
The reductions of the previous section required only the fact that
the coefficients $a^{jj}(x)$ were Lipschitz functions.
The reduction to estimates for angular pieces
depends on the fact that the singularities
of $a^{jj}(x)$, and hence the points
where the $x_2$-derivatives of $p_\l(x,\xi')$ are large, occur
only at $x_2=0$.
Consequently, various error terms that arise in this further reduction
will be highly concentrated at $x_2=0$, which we express through
weighted $L^2$ estimates.

We will take a dyadic decomposition of the $\xi_2$ variable, from
scale $\xi_2\approx\l^{\frac 23}$ to $\xi_2\approx\l$. Thus,
for $1\le j<N_\l=\frac 13\log_2\l$,
let $\beta_j(\xi')=\beta_j(\xi_2,\xi_3)$ denote a smooth cutoff satisfying
$$
\text{supp}(\beta_j)\subset
[2^{-j-2}\l,2^{-j+1}\l]\times [\tfrac 14\l,4\l]\,,
$$
and $\beta_{N_\l}$ supported in
$[-\l^{\frac 23},\l^{\frac 23}]\times[\tfrac 14\l,4\l]\,,$
such that, with $\beta_{-j}(\xi_2,\xi_3)=\beta_j(-\xi_2,\xi_3)$
$$
\sum_{j=1}^{N_\l}\beta_j(\xi')
+\sum_{j=1-N_\l}^{-1}\beta_j(\xi')
=1\quad\text{if}\quad|\xi_2|\le\tfrac 18\l\quad\text{and}\quad \xi_3\in
[\tfrac 12\l,2\l]\,.
$$
Let
$$
u_j(x)=\beta_j(D')u_\l(x)\,.
$$
If we define
$$
\theta_j=2^{-|j|}\,,
$$
then $u_j$ has frequencies localized to
$\xi_2\approx\pm\theta_j\,\xi_3\,,$ or 
$\,|\xi_2|\lesssim \l^{-\frac 13}\xi_3\,$
in case $j=N_\l$.

On the microlocal support of $u_j$, the bicharacteristic equation for
the principal symbol
$\xi_1-p_\l(x,\xi')$ satisfies $\frac{dx_2}{dx_1}\approx\pm\theta_j$,
respectively as $j>0$ or $j<0$.
A bicharacteristic curve passing through the microlocal support
of $u_j$ will satisfy this condition on an interval of $x_1$-length
less than $\eps\theta_j$, if $\eps$ is a small constant.
It is thus natural that we will have good estimates for $u_j$
on slabs of width $\eps\theta_j$ in the $x_1$ variable, and it turns
out this is sufficient to prove Theorem \ref{maintheorem3}.

In proving estimates for $u_j$,
it is convenient to work with the symbol $p_j$
obtained by truncating $p(x,\xi')$ to $x'$-frequencies
less that $c\,\theta_j^{-\frac 12}\l^{\frac 12}\,.$ This finer
truncation than that of $p_\l$ is chosen so that, after
rescaling space by $\theta_j$, the rescaled symbol
$p_j(\theta_j x,\cd)$ will be
$x'$-frequency truncated at $\mu^{\frac 12}$, where 
$\mu=\theta_j\l$ is the frequency scale
of the rescaled solution $u_j(\theta_j x)$. This square root truncation
is consistent with the wave packet techniques we use, and is standard
in the construction of parametrices for rough metrics.

The energy of the induced error term $(P_\l-P_j)u$ will be 
large at $x_2=0$, but decays
away from $x_2=0$ at a rate that is integrable along bicharacteristic curves
that traverse the boundary at angle $\theta_j$. This error term can thus be
considered as a bounded driving force, and we call this term $G_j$ below.

In the next two sections we will establish the following
result.
\begin{theorem}\label{maintheorem4}
Let $S_{j,k}$ denote the slab $x_1\in
[k\eps\theta_j,(k+1)\,\eps\theta_j]\,,$ for
$0\le k\le \eps^{-1}2^{|j|}\,.$

Then, if
$$
D_1u_j-P_j u_j=F_j+G_j\,,
$$
it holds uniformly over $j$ and $k$, and $6\le q\le \infty$, that
\begin{multline*}
\|u_j\|_{L^qL^2(S_{j,k})}\lesssim
\l^{\delta(q)}\theta_j^{\frac 12-\frac 3q}
\Bigl(\;\|u_j\|_{L^\infty L^2(S_{j,k})}+
\|F_j\|_{L^1L^2(S_{j,k})}\\
+\l^{\frac 14}\theta_j^{\frac 14}
\|\langle \l^{\frac 12}\theta_j^{-\frac 12}x_2\rangle^{-1}u_j\|_{L^2(S_{j,k})}+
\l^{-\frac 14}\theta_j^{-\frac 14}
\|\langle\l^{\frac 12}\theta_j^{-\frac 12}x_2\rangle^2 G_j\|_{L^2(S_{j,k})}
\,\Bigr)\,.
\end{multline*}
For $j=N_\l$, it holds that
\begin{equation*}
\|u_j\|_{L^qL^2(S_{j,k})}\lesssim
\l^{\delta(q)}\theta_j^{\frac 12-\frac 3q}
\Bigl(\;\|u_j\|_{L^\infty L^2(S_{j,k})}+
\|F_j+G_j\|_{L^1L^2(S_{j,k})}
\,\Bigr)\,.
\end{equation*}
\end{theorem}

The gain of the factor $\theta_j^{\frac 12-\frac 3q}$ reflects the fact
that, for $q>6$, there is an improvement
in the squarefunction estimates if the solution is localized to a small
conic set in frequency.

The terms $G_j$ arise naturally in both the linearization step
of Lemma \ref{oplemma} and the paradifferential smoothing \eqref{6.2}.
They reflect the fact that the singularities of $d^2 a^{jj}(x)$ are localized
to $x_2=0$. The weighted $L^2$ bound on $u_j$ is a characteristic energy
estimate.

If $\theta_j\approx 1$, then the weighted $L^2$
bound on $G_j$ dominates the $L^1_{x_2}L^2_{x_1,x_3}$ norm of
$G_j$, and exchanging $x_1$ and $x_2$ we could treat $G_j$ and $F_j$
the same. In this case the bound on $u_j$ would be dominated by the
$L^\infty_{x_2}L^2_{x_1,x_3}$ norm.
For small $\theta_j$, however, we cannot use $x_2$ as our ``time''
variable, and we are forced to work with the weighted $L^2$ norms.
These weighted norms can be thought of as an energy
norm along the bicharacteristic flow at angle $\theta_j$. Precisely,
if one replaced $x_2=\theta_j (x_1-c)$ in the weight, then the weighted
$L^2$ norms of $u_j$ and $G_j$ would
behave like the $L^\infty L^2$ and $L^1L^2$ norms
respectively. The crossing point $c$ differs, however, for different
bicharacteristics.

The proof of Theorem \ref{maintheorem4} is contained in sections
\ref{section4} and \ref{section5}. In section \ref{section6}
we establish the appropriate
bounds on the norms occuring on the right side if, as above,
$u_j=\beta_j(D')u_\lambda$, while $F_j$ and $G_j$ are defined in
\eqref{6.1}-\eqref{6.2} below.

To state the bounds required, let $c_{j,k}$ denote the term occuring inside
parentheses on the right hand side of Theorem \ref{maintheorem4}.
In section \ref{section6}, we show that, if $D_1u_\l-P_\l u_\l=F_\l$,
then we have a uniform summability condition
\begin{equation}\label{cijcond}
\sum_j c_{j,k(j)}^2\lesssim \|u_\l\|_{L^\infty L^2(S)}^2+\|F_\l\|_{L^2(S)}^2\,,
\end{equation}
where $k(j)$ denotes any sequence
of values for $k$ such that the slabs $S_{j,k(j)}$ are nested,
in that for $j>0$ we have
$S_{j+1,k(j+1)}\subset S_{j,k(j)}$ (with the analogous condition
for $j<0$.)

In the remainder of this section we show how Theorem \ref{maintheorem3}
follows from Theorem \ref{maintheorem4} together with the bound
\eqref{cijcond}.

We first remark that, if $q$ is a fixed index with $q\ne 8$, 
the bounds of Theorem
\ref{maintheorem3} hold (with constant depending on $q$) under the
weaker assumption that the $c_{j,k}$ are uniformly bounded by the
right side of \eqref{cijcond}. To see this, we sum over the
$2^j\varepsilon^{-1}$ slabs and write
\begin{align*}
\|u_j\|_{L^qL^2(S)}\le\Bigl(\;\sum_{k=1}^{2^j\eps^{-1}}
\|u_j\|_{L^qL^2(S_{j,k})}^q\,\Bigr)^{\frac 1q} &\lesssim
\lambda^{\delta(q)}\theta_j^{\frac12-\frac4q}
\|c_{j,k}\|_{\ell^\infty_j \ell^\infty_k}
\\
&\lesssim \l^{\delta(q)}\theta_j^{\frac 12-\frac 4q}\,
\bigl(\,\|u_\l\|_{L^\infty L^2(S)}+\|F_\l\|_{L^2(S)}\,\bigr)\, .
\end{align*}
The values of $\theta_j=2^{-|j|}$ vary dyadically from $\l^{-\frac
13}$ to $1$. For $q>8$ we can sum over $j$ to obtain
$$
\|u_\l\|_{L^qL^2(S)}\lesssim\l^{\delta(q)}
\bigl(\,\|u_\l\|_{L^\infty L^2(S)}+\|F_\l\|_{L^2(S)}\,\bigr)\,,
$$
and for $6\le q<8$ the sum yields
$$
\|u_\l\|_{L^qL^2(S)}\lesssim\l^{\delta(q)-\frac 13(\frac 12-\frac 4q)}
\bigl(\,\|u_\l\|_{L^\infty L^2(S)}+\|F_\l\|_{L^2(S)}\,\bigr)\,.
$$
The above exponent of $\l$ equals $\gamma(q)$, yielding the desired bound.
The geometric sum, however, increases as $q\rightarrow 8$, and yields
a logarithmic loss in $\l$ at $q=8$.

To obtain the bound at $q=8$, and hence uniform bounds over $q$ in
Theorem \ref{maintheorem3}, we use the following
worst-case branching argument. We consider terms with $j>0$ here,
the negative terms being controlled by the same argument.

Let $S_{1,k(1)}$ denote the slab at scale $\eps\,2^{-1}$
that maximizes $\|u_\l\|_{L^8L^2(S_{1,k})}\,.$
Since the decomposition of $u_\l$ into $u_j$ is a Littlewood-Paley
decomposition in the $\xi_2$ variable, we have
$$
\|u_\l\|_{L^8L^2(S_{1,k(1)})}^2\lesssim
\Bigl\|\Bigl(\;\sum_{j=1}^{N_\l}|u_j|^2\Bigr)^{\frac 12}
\Bigr\|_{L^8L^2(S_{1,k(1)})}^2\,.
$$
By the Minkowski inequality,
\begin{align*}
\Bigl\|\Bigl(\;\sum_{j=1}^{N_\l}|u_j|^2\Bigr)^{\frac 12}
\Bigr\|_{L^8L^2(S_{1,k(1)})}^2
&\le \|u_1\|_{L^8L^2(S_{1,k(1)})}^2+
\Bigl(\;\sum_{S_{2,k}\subset S_{1,k(1)}}\!\!\!
\Bigl\|\Bigl(\;\sum_{j=2}^{N_\l}|u_j|^2\Bigr)^{\frac 12}
\Bigr\|_{L^8L^2(S_{2,k})}^8
\;\Bigr)^{\frac 28}\\
&\le  \|u_1\|_{L^8L^2(S_{1,k(1)})}^2+
2^{\frac 28}\,
\Bigl\|\Bigl(\;\sum_{j=2}^{N_\l}|u_j|^2\Bigr)^{\frac 12}
\Bigr\|_{L^8L^2(S_{2,k(2)})}^2
\end{align*}
where $k(2)$ is chosen to maximize
$\bigl\|\bigl(\sum_{j=2}^\infty|u_j|^2\bigr)^{\frac 12}\|_{L^8(S_{2,k})}$
among the two slabs $S_{2,k}$ contained in $S_{1,k(1)}$.
Repeating this procedure yields a nested sequence such that
\begin{align*}
\eps^{\frac 14}\|u_\l\|_{L^8L^2(S)}^2&\le
\|u_1\|_{L^8L^2(S_{1,k(1)})}^2+
2^{\frac 28}\|u_2\|_{L^8L^2(S_{2,k(2)})}^2+
2^{\frac 48}\|u_3\|_{L^8L^2(S_{3,k(3)})}^2+\,\cdots\\
&\lesssim \l^{2\delta(8)}\bigl(\,
 c_{1,k(1)}^2+c_{2,k(2)}^2+c_{3,k(3)}^2+\,\cdots\,\bigr)
\end{align*}
where the last holds by Theorem \ref{maintheorem4}
since $\theta_j^{\frac 12-\frac 38}=2^{-\frac j8}\,.$
The case $q=8$ of Theorem \ref{maintheorem3}
follows by \eqref{cijcond}.\qed

\newsection{The Wave Packet Transform}\label{section4}
The purpose of this section and the next is to establish
Theorem \ref{maintheorem4}. We assume for these two sections that
we have fixed $\l$ and $\theta_j$, and consider $j>0$ so that
$\xi_2>0$ (except for the term $j=N_\l$, where $|\xi_2|\le\l^{\frac 23}\,.$)

We will rescale space by $\theta_j$. Thus, we work with the function
$$
u(x)=u_j(\theta_j x)\,,
$$
which for $j\ne N_\l$ is supported in the set
$$
\xi_2\in \bigl[\tfrac 14\theta_j\mu,2\theta_j\mu\bigr]\,,
\qquad\xi_3\in\bigl[\tfrac 14\mu,4\mu\bigr]\,,
$$
where
$$
\mu=\theta_j\l
$$
is now the frequency scale for $u(x)$. For $j=N_\l$, we have
$|\xi_2|\le\mu^{\frac 12}$, and $\theta_{N_\l}=\mu^{-\frac 12}\,.$

Let $q(x,\xi')$ denote the rescaled symbol
$$
q(x,\xi')=\theta_j p_j(\theta_j x,\theta_j^{-1}\xi')\,,
$$
which is truncated to $x'$-frequencies less than $c\mu^{\frac 12}$.
For $|\xi'|\approx\mu$, the symbol $q$ satisfies the estimates
\begin{equation}\label{qsymbolcond}
\bigl|\partial_x^\beta\partial^\alpha_{\xi'} q(x,\xi')\bigr|
\lesssim\begin{cases}
\mu^{1-|\alpha|}\,,\qquad |\beta|=0\,,\\
c_0\,
\bigl(\,1+\mu^{\frac 12(|\beta|-1)}\theta_j\,\langle
\mu^{\frac 12}x_2\rangle^{-N}\,\bigr)
\,\mu^{1-|\alpha|}\,,\quad|\beta|\ge 1\,.
\end{cases}
\end{equation}
This follows from \eqref{qsymbolcond'}.

In the remainder of this section and the next, we
will drop the index $j$. The quantities $\theta$ and $\mu$ are the
two relevant parameters for our purposes.
After rescaling the estimates of Theorem \ref{maintheorem4},
and translating $S_{j,k}$ in $x_1$ to $x_1=0$,
we are reduced to establishing the following. Here, $S$ denotes the
$(x_1,x')$ slab $[0,\eps]\times\R^2$.

\begin{theorem}\label{maintheorem5}
Suppose that $\widehat u(\xi)$ is supported in the set
$$
\xi_2\in \bigl[\tfrac 14\theta\mu,2\theta\mu\bigr]\,,\qquad
\xi_3\in \bigl[\tfrac 14\mu,4\mu\bigr]\,,
$$
respectively $|\xi_2|\le\mu^{\frac 12}$ in case
$\theta=\mu^{-\frac 12}\,.$ Suppose that $u$ satisfies
$$
D_1 u-q(x,D')u=F+G
$$
on the slab $S$, where $q$ satisfies \eqref{qsymbolcond}, and is
truncated to $x'$-frequencies less than $c\mu^{\frac 12}$. Then
the following bounds hold,
uniformly over $\theta$ and $\mu$, and $6\le q\le\infty$,
\begin{multline*}
\|u\|_{L^qL^2(S)}\lesssim
\mu^{\delta(q)}\theta^{\frac 12-\frac 3q}
\Bigl(\;\|u\|_{L^\infty L^2(S)}+
\|F\|_{L^1L^2(S)}\\
+\mu^{\frac 14}\theta^{\frac 12}
\|\langle \mu^{\frac 12}x_2\rangle^{-1}u\|_{L^2(S)}+
\mu^{-\frac 14}\theta_j^{-\frac 12}
\|\langle\mu^{\frac 12}x_2\rangle^2 G\|_{L^2(S)}
\,\Bigr)\,,
\end{multline*}
and for $\theta=\mu^{-\frac 12}$
\begin{equation*}
\|u\|_{L^qL^2(S)}\lesssim
\mu^{\delta(q)}\theta^{\frac 12-\frac 3q}
\Bigl(\;\|u\|_{L^\infty L^2(S)}+
\|F+G\|_{L^1L^2(S)}
\,\Bigr)\,.
\end{equation*}
\end{theorem}

We introduce at this point the phase-space transform which will be used
to establish Theorem \ref{maintheorem5}. This transform is essentially
the C\'ordoba-Fefferman wave packet transform \cite{CF}. The precise
form here is a simple modification of the FBI transform
used by Tataru in \cite{Tat1} and \cite{Tat3}
to establish Strichartz estimates for low regularity metrics;
the difference is that in our applications
we use a Schwartz function with compactly supported Fourier transform,
instead of the Gaussian, as the fundamental
wave packet. This is useful in that it strictly localizes the frequency
support of the transformed functions.
Our transform will act on the $x'=(x_2,x_3)$ variables.

We use the notion of previous sections: $x=(x_1,x_2,x_3)=(x_1,x')$,
where $x_3$ denotes the variable $t$.

Fix a real, radial Schwartz function $g(z')\in S(\R^2)$, with
$\|g\|_{L^2(\R^2)}=(2\pi)^{-1}$,
and assume that its Fourier transform
$\widehat{g}(\zeta')$
is supported in the ball $\{|\zeta'|\le c\}\,.$
For $\mu\ge 1$, we define
$T_\mu\,:\, S'(\R^2)\rightarrow C^\infty(\R^4)$ by the rule
$$
\bigl(T_\mu f\bigr)(x',\xi')=\mu^{\frac 12}\,
\int e^{-i\langle \xi',y'-x'\rangle}\,g\bigl(\mu^{\frac 12}(y'-x')\bigr)\,f(y')
\,dy'\,.
$$
A simple calculation shows that
$$
f(y')=\mu^{\frac 12}\,
\int e^{i\langle \xi',y'-x'\rangle}\,g\bigl(\mu^{\frac 12}(y'-x')\bigr)\,
\bigl(T_\mu f\bigr)(x',\xi')\,dx'\,d\xi'\,,
$$
so that $T_\mu^*T_\mu=I\,.$ In particular,
\begin{equation}\label{fest}
\|T_\mu f\|_{L^2(\R^4)}=\|f\|_{L^2(\R^2)}\,.
\end{equation}
It will be useful to note that this holds in a more general setting.
\begin{lemma}\label{wavetype}
Suppose that $g_{x',\xi'}(y')$ is a family of Schwartz functions on $\R^2$,
depending on the parameters $x'$ and $\xi'$, with uniform bounds
over $x'$ and $\xi'$ on each Schwartz norm of $g$. Then the operator
$$
\bigl(T_\mu f\bigr)(x',\xi')=\mu^{\frac 12}\,
\int e^{-i\langle \xi',y'-x'\rangle}\,g_{x',\xi'}
\bigl(\mu^{\frac 12}(y'-x')\bigr)\,f(y')\,dy'
$$
satisfies the bound
$$
\|T_\mu f\|_{L^2(\R^4)}\lesssim \|f\|_{L^2(\R^2)}\,.
$$
\end{lemma}
\demo $T_\mu$ is bounded if and only if $T^*_\mu$ is bounded.
Since $\|T^*_\mu F\|_{L^2}^2\le \|T_\mu T^*_\mu F\|_{L^2}\|F\|_{L^2}$,
it suffices to see that $T_\mu T^*_\mu$ is bounded on
$L^2(dy'd\xi')$.

The operator $T_\mu T_\mu^*$ is an integral operator with kernel
$$
K(y',\eta';x',\xi')=
\mu\,e^{i\langle\eta',y'\rangle-i\langle\xi',x'\rangle}
\int e^{i\langle\xi'-\eta',z'\rangle}
g_{y',\eta'}\bigl(\mu^{\frac 12}(z'-y')\bigr)\,
\overline{g_{x',\xi'}\bigl(\mu^{\frac 12}(z'-x')\bigr)}
\,dz'\,.
$$
A simple integration by parts argument shows that
$$
\bigl|K(y',\eta';x',\xi')\bigr|\lesssim
\bigl(\,1+\mu^{-\frac 12}|\eta'-\xi'|+\mu^{\frac 12}|y'-x'|\,
\bigr)^{-N}\,,
$$
with constants depending only on uniform bounds for a finite
collection of seminorms of $g_{x',\xi'}$ depending on $N$. The
$L^2(\R^4)$ boundedness of $T_\mu T^*_\mu$ then follows by Schur's
Lemma. \findemo

A corollary of this lemma is that, for $N$ positive or negative,
\begin{equation}\label{weightest}
\|\langle \mu^{\frac 12}x_2\rangle^N T_\mu f\|_{L^2(\R^4)}\lesssim
\|\langle \mu^{\frac 12}x_2\rangle^N f\|_{L^2(\R^2)}\,,
\end{equation}
by considering $g_{x'}(y)=\langle \mu^{\frac 12}x_2\rangle^N
\langle \mu^{\frac 12}x_2-y_2\rangle^{-N}g(y')\,.$

The next two lemmas relate the conjugation of $q(x,D')$, by the wave
packet transform, to the Hamiltonian flow under $q$. 
These results are analogous to Theorem 1 of \cite{Tat1}, in that the 
error term is one order better where the metric has two bounded
derivatives. In our case the second derivatives are large along
the boundary $x_2=0$, which leads to larger errors there. A key fact for
our paper is that the errors are suitably integrable along the 
Hamiltonian flow of $q$.

\begin{lemma}\label{pest}
Let $q(x,\xi')$ satisfy the estimates \eqref{qsymbolcond}.
Suppose that $|\xi'|\approx\mu\,.$ Then, if $q(y,\Dy')^*$ acts on the
$y'$ variable, and $y_1=x_1$, we can write
\begin{multline*}
\Bigl(\,q(y,\Dy')^*-id_{\xi'} q(x,\xi')\cdot d_{x'}+id_{x'} q(x,\xi')\cdot
d_{\xi'}\,\Bigr)
\Bigl[e^{i\langle \xi',y'-x'\rangle}\,g\bigl(\mu^{\frac 12}(y'-x')\bigr)\Bigr]
\\
=
e^{i\langle \xi',y'-x'\rangle}\,g_{x,\xi'}\bigl(\mu^{\frac 12}(y'-x')\bigr)
\end{multline*}
where $g_{x,\xi'}(\cd)$ denotes a family of Schwartz functions on $\R^2$
depending on the parameters $x$ and $\xi'$, each of which has
Fourier transform supported in the ball of radius $2c$.
If $\|\cd\|$ denotes any of the Schwartz seminorms,
we have
$$
\|g_{x,\xi'}\|\lesssim
1+c_0\,\mu^{\frac 12}\theta\langle\mu^{\frac 12}x_2\rangle^{-3}
\,,
$$
where $c_0$ is the small constant of \eqref{cond3}.
\end{lemma}
\demo
Letting $\mathfrak F$ denote the Fourier transform with respect to $y'$,
we write
\begin{multline*}
{\mathfrak F}\circ
\Bigl(q(y,\Dy')^*-id_{\xi'} q(x,\xi')\cdot d_{x'}
+id_{x'} q(x,\xi')\cdot d_{\xi'}\,
\Bigr)
\Bigl[e^{i\langle \xi',y'-x'\rangle}
\,g\bigl(\mu^{\frac 12}(y'-x')\bigr)\,\Bigr](\eta')\\
=
e^{-i\langle\eta',x'\rangle}\,\mu^{-1}
\,\widehat{g_{x,\xi'}}\bigl(\mu^{-\frac 12}(\eta'-\xi')\bigr)\,,
\end{multline*}
where $\widehat{g_{x,\xi'}}(\eta')$ is equal to
\begin{multline*}
\int e^{-i\langle\eta',y'\rangle}
\Bigl[\,q(x+\mu^{-\frac 12}y',\xi'
+\mu^{\frac 12}\eta')-q(x,\xi')
-d_{x',\xi'} q(x,\xi')\cdot\bigl(\mu^{-\frac 12}y',\mu^{\frac 12}\eta'\bigr)
\,\Bigr]\,g(y')\,dy'\\
=
\int_0^1(1-\sigma)\,\left[\int  e^{-i\langle\eta',y'\rangle}
\partial_\sigma^{\,2}\Bigl(\,
q\bigl(x+\sigma\mu^{-\frac 12}y',\xi'+\sigma\mu^{\frac 12}\eta'\bigr)\Bigr)
\,g(y')
\,dy'\,\right]d\sigma\,.
\end{multline*}
The spectral restriction on $q$ and $g$ imply that this vanishes for
$|\eta'|\ge 2c\,.$ Consequently, it suffices to establish
$C^\infty$ bounds in $\eta'$ for the term in brackets, uniformly over
$\sigma\in [0,1]$ and $|\eta'|\le 2c\,.$
Since the effect of differentiating the integrand with respect to $\eta'$
is innocuous, as the rapid decrease in $g(y')$ counters any polynomial
in $y'$,  we content ourselves with establishing uniform pointwise
bounds on the term in brackets. Note that
$|\xi'+\sigma\mu^{\frac 12}\eta'|\approx \mu\,.$

The effect of $\partial_\sigma^2$ is to bring out factors of
$\mu^{\pm\frac 12}$, and to differentiate $q$ twice.
If $q$ is differentiated at most once in $x'$, then the bounds
$$
|\partial_{x'}\partial_{\xi'} q(x,\xi')|\lesssim 1\,,\qquad
|\partial_{\xi'}^2 q(x,\xi')|\lesssim \mu^{-1}\,,\quad
{\rm for}\quad
|\xi'|\approx \mu\,,
$$
yield bounds of size $1$ on the term. If $q$ is differentiated twice in
$x'$, then by \eqref{qsymbolcond} we have the bounds,
for $|\xi'|\approx \mu\,,$
\begin{align*}
\mu^{-1}|\partial_{x'}^2 q(x+\sigma\mu^{-\frac 12}y',\xi')|\lesssim&\,
c_0+c_0\,\mu^{\frac 12}\theta\langle\mu^{\frac 12}x_2+\sigma y_2\rangle^{-3}\\
\lesssim&\,
1+c_0\,\mu^{\frac 12}\theta\langle\mu^{\frac 12}x_2\rangle^{-3}
\langle y_2\rangle^3\,.
\end{align*}
The rapid decrease of $g(y')$ absorbs the term
$\langle y_2\rangle^3$, leading
to the desired bounds.
\findemo

We now take the wave packet transform of the solution $u(x)$ with
respect to the $x'$ variables, and introduce the notation
$\tilde{u}(x,\xi')=(T_\mu u)(x,\xi')\,.$ The functions
$\tilde{F}(x,\xi')$ and $\tilde{G}(x,\xi')$ in the next lemma,
though, include terms in addition to the transforms of $F$ and $G$
of Theorem \ref{maintheorem5}. Let $\tilde S$ denote the
$(x_1,x',\xi')$ slab $[0,\eps]\times\R^4=S\times\R^2_{\xi'}$.
\begin{lemma}\label{oplemma}
Under the above conditions, we may write
$$
\Bigl(d_1-d_{\xi'} q(x,\xi')\cdot d_{x'}+d_{x'} q(x,\xi')\cdot d_{\xi'}\Bigr)
\tilde{u}(x,\xi')
=\tilde{F}(x,\xi')+\tilde{G}(x,\xi')\,,
$$
where
\begin{multline}\label{FGest}
\|\tilde{F}\|_{L^1L^2(\tilde S)}+
\mu^{-\frac 14}\theta^{-\frac 12}
\|\langle\mu^{\frac 12}x_2\rangle^2 \tilde{G}
\|_{L^2(\tilde S)}
\\
\lesssim\,
\|u\|_{L^\infty L^2(S)}+
\|F\|_{L^1 L^2(S)}
+\mu^{\frac 14}\theta^{\frac 12}
\|\langle \mu^{\frac 12}x_2\rangle^{-1}u\|_{L^2(S)}+
\mu^{-\frac 14}\theta^{-\frac 12}
\|\langle\mu^{\frac 12}x_2\rangle^2 G\|_{L^2(S)}\,.
\end{multline}
Furthermore, $\tilde F$ and $\tilde G$ are supported in a set where
$\xi_2\approx\theta\mu\,,$
$\,\xi_3\approx\mu\,.$

In case $\theta=\mu^{-\frac 12}$, then
\begin{equation}\label{FGest2}
\|\tilde{F}+\tilde{G}\|_{L^1L^2(\tilde S)}
\lesssim\,
\|u\|_{L^\infty L^2(S)}+
\|F+G\|_{L^1L^2(S)}\,,
\end{equation}
and $\tilde{F}+\tilde{G}$ is supported where
$|\xi_2|\lesssim \mu^{\frac 12}$ and $\xi_3\approx\mu$.
\end{lemma}
\demo
Applying $T_\mu$ to the equation
$D_1u=F+G+q(x,D')u$ yields
\begin{multline*}
d_1 \tilde{u}(x,\xi')=i\bigl(T_\mu F\bigr)(x,\xi')+
i\bigl(T_\mu G\bigr)(x,\xi')
\\+i\,
\mu^{\frac 12}\!
\int \;\overline{q(x_1,y',\Dy')^*
\Bigl[\,e^{i\langle \xi',y'-x'\rangle}
\,g\bigl(\mu^{\frac 12}(y'-x')\bigr)\,\Bigr]}
\,u(x_1,y')\,dy'\,.
\end{multline*}
The terms $T_\mu F$ and $T_\mu G$ satisfy the bounds required of $\tilde F$
and $\tilde G$ respectively, the latter by the estimate \eqref{weightest}
in the case of \eqref{FGest}.
By Lemma \ref{pest}, we can write the last term as
\begin{multline*}
\Bigl(d_{\xi'} q(x,\xi')\cdot d_{x'}-d_{x'} q(x,\xi')\cdot d_{\xi'}\Bigr)
\tilde{u}(x,\xi')
\\
+\mu^{\frac 12}\,\int
e^{-i\langle \xi',y'-x'\rangle}\,g_{x,\xi'}\bigl(\mu^{\frac 12}(y'-x')\bigr)
\,u(x_1,y')\,dy'\,.
\end{multline*}
For $x_2$ such that
$\mu^{\frac 12}\theta\langle\mu^{\frac 12} x_2\rangle^{-3}\le 1$, the
latter term is absorbed into $\tilde F$ by Lemmas \ref{wavetype}
and \ref{pest}.
For $x_2$ such that
$\mu^{\frac 12}\theta\langle\mu^{\frac 12} x_2\rangle^{-3}\ge 1$, the
term can be absorbed into $\tilde G$, by \eqref{weightest} and
Lemma \ref{pest}. Here we use the simple fact that \eqref{weightest} holds
for operators of the type in Lemma \ref{wavetype}. Note that if
$\theta=\mu^{-\frac 12}$ the entire term can be absorbed into $\tilde F$.

The support condition on $\tilde F$ and $\tilde G$ follows from the
support condition on $\widehat{u}$, and the fact that $g_{x,\xi'}$
has Fourier
transform supported in the ball of radius $2c$. Alternatively,
we may multiply both sides of the equation defining $\tilde F+\tilde G$
by a cutoff supported in the set $\xi_3\approx\mu$, $\xi_2\approx\theta\mu$
(respectively $|\xi_2|\lesssim\mu^{\frac 12}$),
which equals 1 on the support of $\tilde u$.
\findemo

Let $\Theta_{s,r}$ denote the canonical transform on
$\R^4_{x',\xi'}=T^*(\R^2_{x'})$ generated by the Hamiltonian flow
of $q$. Thus, $\Theta_{s,r}(x',\xi')=\gamma(s)$, where $\gamma$ is
the integral curve of the vector field
$$
d_1-d_{\xi'}q(x,\xi')\cdot d_{x'}+d_{x'}q(x,\xi')\cdot d_{\xi'}
$$
with $\gamma(r)=(x',\xi')$. Note that $\Theta_{s,r}$ is symplectic,
thus preserves the measure $dx'\,d\xi'$, hence induces a unitary
mapping on $L^2(\R^4)$. Furthermore, $\Theta_{s,r}$ maps a set of the
form $\xi_2\approx \theta\xi_3$ to a set of similar form, provided
$|s-r|\le 1$. This follows since $|d_{x'}q(x,\xi')|\le
c\,\theta\,|\xi'|$ for $c$ a small constant.

We can now write
\begin{equation}\label{udecomp}
\tilde{u}(x,\xi')=\tilde{u}(0,\Theta_{0,x_1}(x',\xi'))+ \int_0^{x_1}
\tilde{F}(s,\Theta_{s,x_1}(x',\xi'))\,ds +\int_0^{x_1}
\tilde{G}(s,\Theta_{s,x_1}(x',\xi'))\,ds\,.
\end{equation}
By the preceding comments, for each $s$ the integrands are
supported in the flowout under $\Theta_{s,0}$ of the set
$\xi_3\approx\mu$, $\xi_2\approx \theta\mu$ (respectively
$|\xi_2|\lesssim\mu^{\frac 12}$).

Writing $u=T_\mu^*\tilde{u}$ shows that $u(x)$ can be written as a
superposition of functions, each of which is the restriction to $x_1>s$
of the image under
$T_\mu^*$ of a function invariant under the Hamiltonian flow of
$q$. However, in view of the bounds \eqref{FGest},
the term $\tilde G$ may have large $L^1L^2$ norm if $\theta$ is small,
in contrast to the setting of \cite{Tat3}.
As a result, one cannot directly apply \eqref{udecomp} to reduce
matters to considering estimates for such flow-invariant
functions.
Nevertheless, we can use arguments from
Koch-Tataru \cite{KT} together with \eqref{FGest}
to see that we may indeed reduce
consideration to the case that $\tilde{u}$ is invariant under the
flow of $q$. Roughly, the $V^2_q$ space of \cite{KT} permits
us to use the weaker condition of
integrability of $\tilde G$ along the flow lines of $q$.
We show here that Theorem
\ref{maintheorem5} is a consequence of the following theorem,
which will be proven in the next section.

\begin{theorem}\label{maintheorem6}
Suppose that $f\in L^2(\R^4)$ is supported in
a set of the form $\xi_3\approx\mu\,,$ $\xi_2\approx\theta\mu\,,$
or a set of the form
$\xi_3\approx\mu\,,$ $|\xi_2|\lesssim\mu^{\frac 12}$ in case
$\theta=\mu^{-\frac 12}$.

Then, if
$u=T^*_\mu\bigl[f\bigl(\Theta_{0,x_1}(x',\xi')\bigr)\bigr]\,,$ we
have for $q\ge 6$
$$
\|u\|_{L^qL^2(S)}
\lesssim\mu^{\delta(q)}\theta^{\frac 12-\frac 3q}
\,\|f\|_{L^2(\R^4)}\,.
$$
\end{theorem}

In the remainder of this section we demonstrate the
reduction of Theorem \ref{maintheorem5} to Theorem \ref{maintheorem6}.
In the case of $\theta=\mu^{-\frac 12}$, it is a simple consequence of
\eqref{FGest2} and \eqref{udecomp}.

For general $\theta$ this reduction requires the introduction of
the space $V^2_q$ of functions on $\tilde S$ with bounded
2-variation along the Hamiltonian flow of $q$. Recall that
$\Theta_{r,s}$ preserves the measure $dx'\,d\xi'$. Then,
following Koch-Tataru \cite{KT} we define
$$
\|\tilde{u}\|^2_{V^2_q}=\|\tilde{u}(0,\cd)\|^2_{L^2(\R^4)}+
\sup_{P}
\sum_{j\ge 1}\,
\|\tilde{u}(s_j,\cd)-\tilde{u}(s_{j-1},\Theta_{s_{j-1},s_j}(\cd))
\|^2_{L^2(\R^4)}\,,
$$
where $P$ denotes the family of finite partitions
$\{0=s_0<s_1<\ldots<s_m=\eps\}$ of $[0,\eps]$.

By Lemma 6.4 of \cite{KT}, if $\|\tilde{u}\|_{V^2_q}<\infty$, we may decompose
$$
\tilde{u}=\sum_{k=1}^\infty c_k\,\tilde{u}_k\,,
\qquad\text{with}\quad\sum_{k=1}^\infty |c_k|
\le \|\tilde{u}\|_{V^2_q}\,,
$$
where each function $\tilde{u}_k$ is an atom, in the sense that
for some partition $\{s_j\}$ in $P$
$$
\tilde{u}_k(x,\xi')= \sum_{j=1}^{m}
1_{[s_{j-1},s_j)}(x_1)f_j(\Theta_{0,x_1}(x',\xi'))\,,
$$
where, for each $q>2$, it holds that
$$
\Bigl(\; \sum_{j=1}^m \|f_j\|_{L^2(\R^4)}^q\Bigr)^{\frac 1q}
\le C_q\,.
$$
Note that one may bound $C_q\le C_6$ for $q\ge 6$, so we may take
$C_q$ uniformly bounded, since we work with $q\ge 6$.

We also note that each $f_j$ arising in the atomic decomposition of $\tilde u$
will be supported in the region
$\xi_3\approx\mu\,,$ $\xi_2\approx\theta\mu\,.$ This follows from
the inductive construction of $f_j$ in \cite{KT},
together with the comments surrounding \eqref{udecomp}.

Consider $u_k=T^*_\mu\tilde{u}_k\,.$
Then, assuming Theorem \ref{maintheorem6}, for $q\ge 6$ we may bound
$$
\|u_k\|_{L^qL^2(S)}\le \Bigl(\;\sum_{j=1}^m
\|T^*_\mu\bigl[f_j\bigl(\Theta_{0,x_1}(\cd)\bigr)\bigr]
\|_{L^qL^2(S)}^q\Bigr)^{\frac 1q}\lesssim \Bigl(\;\sum_{j=1}^m
\|f_j\|_{L^2(\R^4)}^q\Bigr)^{\frac 1q} \lesssim 1\,.
$$
Summing over $k$ yields $\|u\|_{L^qL^2(S)}
\lesssim\|\tilde u\|_{V^2_q}\,.$
It thus remains to demonstrate that
\begin{equation}\label{uv2bound}
\|\tilde u\|_{V^2_q}\lesssim
\|\tilde u(0,\cd)\|_{L^2(\R^4)}+
\|\tilde{F}\|_{L^1L^2(\tilde S)}
+\mu^{-\frac 14}\theta^{-\frac 12}\|\langle\mu^{\frac 12}x_2\rangle^2 \tilde{G}
\|_{L^2(\tilde S)}
\,,
\end{equation}
since by Lemma \ref{oplemma} and boundedness of $T_\mu$ the right
hand side here
is dominated by the right hand side in Theorem \ref{maintheorem5}.

We use the decomposition \eqref{udecomp}, and note that the $V^2_q$
norm of the first two
terms on the right hand side of \eqref{udecomp} are easily bounded
by the first two terms on the right hand side of \eqref{uv2bound},
the latter since
\begin{align*}
\sum_j\Bigl\|\int_0^{s_j}\tilde F(s,\Theta_{s,s_j}(x',\xi'))&\,ds
-\int_0^{s_{j-1}}\tilde
F(s,\Theta_{s,s_j}(x',\xi'))\,ds\Bigr\|_{L^2({\mathbb R}^4)}^2
\\
&\le \sum_j\left(\int_{s_{j-1}}^{s_j}\|\tilde
F(s,\Theta_{s,s_j}(x',\xi'))\|_{L^2({\mathbb R}^4)}\, ds\right)^2
\\
&
=\sum_j \left(\int_{s_{j-1}}^{s_j}\|\tilde
F(s,\cdot)\|_{L^2({\mathbb R}^4)} \, ds\right)^2\lesssim \|\tilde
F\|^2_{L^1L^2(\tilde S)},
\end{align*}
using the invariance of $dx'd\xi'$ under $\Theta$ in the second equality.

We thus reduce to the case that $\tilde F=0$ and $\tilde u(0,x',\xi')=0$,
and hence by \eqref{udecomp} that
$$
\tilde{u}(x,\xi')= \int_0^{x_1}
\tilde{G}(s,\Theta_{s,x_1}(x',\xi'))\,ds\,.
$$
Note that, by the group property of $\Theta$, we have
\begin{equation}\label{u2var}
\|\tilde u(s_j,\cd)-\tilde u(s_{j-1},\Theta_{s_{j-1},s_j}(\cd))
\|^2_{L^2(\R^4)}= \Bigl\|\;\int_{s_{j-1}}^{s_j}\tilde
G(s,\Theta_{s,s_j}(\cd))\,ds\, \Bigr\|^2_{L^2(\R^4)}\,.
\end{equation}
Given a partition $\{0=s_0<s_1<\ldots<s_m=\eps\}$, we first
consider the sum of the quantity \eqref{u2var} over those indices
$j$ for which $|s_j-s_{j-1}|\le \mu^{-\frac 12}\theta^{-1}$. By
the Schwarz inequality we may bound the sum by
$$
\sum_j \theta^{-1}\mu^{-\frac 12} \|\tilde
G(s,\Theta_{s,s_j}(x',\xi')) \|_{L^2([s_{j-1},s_j]\times\R^4)}^2 \le
\mu^{-\frac 12}\theta^{-1}\|\tilde G\|_{L^2(\tilde S)}^2\,.
$$
%\,,
%$$
%where in the last step we used invariance of $dx'\,d\xi'$ under
%$\Theta$.

Next, consider an index $j$ for which
$|s_j-s_{j-1}|>\mu^{-\frac 12}\theta^{-1}\,.$ We split the
interval $[s_j,s_{j-1}]$ into a union of intervals $I_k$
for which $\frac 12|I_k|\le\mu^{-\frac 12}\theta^{-1}\le |I_k|\,.$
We claim that we may bound
\begin{equation}\label{ortho}
\Bigl\|\;\int_{s_{j-1}}^{s_j}\tilde G(s,\Theta_{s,s_j}(\cd))\,ds\,
\Bigr\|^2_{L^2(\R^4)}\lesssim \sum_k\; \Bigl\|\;\langle\mu^{\frac
12}x_2\rangle^2 \int_{I_k}\tilde G (s,\Theta_{s,s_k}(\cd))\,ds\,
\Bigr\|^2_{L^2(\R^4)}\,,
\end{equation}
where $s_k$ denotes the right endpoint of $I_k$. Given
\eqref{ortho}, we may apply the Schwarz inequality as before
(together with the fact that the weight $\langle\mu^{\frac
12}x_2\rangle^2$ is essentially preserved by $\Theta_{s,s_k}$, since
$\frac{dx_2}{dx_1}\approx\theta$ on the domain of integration and
$|s-s_k|\le \mu^{-\frac 12}\theta^{-1}$) to bound the sum over $k$
and then $j$ by the right hand side of \eqref{uv2bound}.

To prove \eqref{ortho}, we write
$$
\int_{s_{j-1}}^{s_j}\tilde G(s,\Theta_{s,s_j}(x',\xi'))\,ds= \sum_k
\tilde v_k(\Theta_{s_k,s_j}(x',\xi'))\,,
$$
with
$$
\tilde v_k(x',\xi')=\int_{I_k}\tilde G(s,\Theta_{s,s_k}(x',\xi'))\,ds
\,.
$$
Then \eqref{ortho} will follow by showing that
\begin{multline*}
\Bigl|\;\int \tilde v_k(\Theta_{s_k,s_j}(x',\xi'))\, \overline{\tilde
v_{k'}(\Theta_{s_{k'},s_j}(x',\xi'))}\,dx'\,d\xi'\,
\Bigr|\\
=\Bigl|\;\int \tilde v_k(\Theta_{s_k,s_{k'}}(x',\xi'))\,
\overline{\tilde v_{k'}(x',\xi')}\,dx'\,d\xi'\,
\Bigr|\\
\lesssim
|k-k'|^{-2}\,
\|\langle\mu^{\frac 12}x_2\rangle^2\tilde v_k\|_{L^2(\R^4)}
\|\langle\mu^{\frac 12}x_2\rangle^2\tilde v_{k'}\|_{L^2(\R^4)}\,.
\end{multline*}
This, in turn, is a simple consequence of the fact that
$\frac{dx_2}{dx_1}\approx \theta$ on the domain of integration,
and hence, letting $x_2$ denote the $x_2$-coordinate function,
$$
|\,x_2(\Theta_{s_k,s_{k'}}(x',\xi'))-x_2|\approx
\theta\,|s_k-s_{k'}|\approx \mu^{-\frac 12}|k-k'|\,.
$$
Consequently,
$$
\langle\mu^{\frac 12}x_2(\Theta_{s_k,s_{k'}}(x',\xi'))\rangle^{-2}
\langle\mu^{\frac 12}x_2\rangle^{-2}\lesssim |k-k'|^{-2}\,.\qed
$$

\newsection{Homogeneous Estimates}\label{section5}
In this section we prove Theorem \ref{maintheorem6}.
For notational convenience,
the variables $z=(z_2,z_3)$ and $\zeta=(\zeta_2,\zeta_3)$
will be used as dummy
variables in the role of $x'$ and $\xi'$, as will $w$ and $\eta$.
We also use real variables $r,s,t$ as dummy variables in the role
of $x_1$ and $y_1$.
For $f\in L^2(dx'\,d\xi')$, define $W\!f$ by the rule
$$
W\!f(x)=T_\mu^*\bigl(f\circ \Theta_{0,x_1}\bigr)(x')\,.
$$
Let $\beta_\theta(\xi')$ be a cutoff to the region $\xi_3\approx\mu$,
$\xi_2\approx \theta\mu$, (respectively $|\xi_2|\le\mu^{\frac 12}$
in case $\theta=\mu^{-\frac 12}$).
Then Theorem \ref{maintheorem6} is equivalent to establishing, for $q\ge 6$,
the bound
$$
\|\beta_\theta(D')W\!f\|_{L^qL^2(S)}\lesssim
\mu^{\delta(q)}\theta^{\frac 12-\frac 3q}\|f\|_{L^2(\R^4)}\,,
$$
which is equivalent to the bound
\begin{equation}\label{Fbound}
\|\beta_\theta(D')WW^*\beta_\theta(D')F\|_{L^qL^2(S)}\lesssim
\mu^{2\delta(q)}\theta^{1-\frac 6q}
\|F\|_{L^{q'}L^2(S)}\,.
\end{equation}
The operator $WW^*$ takes the form
$$
\bigl(WW^*\!F\bigr)(x)= \int_{0}^\eps T_\mu^*\bigl[\bigl(T_\mu
F\bigr)(s,\cd)\circ\Theta_{s,x_1}\bigr](x')\,ds\,.
$$
If applied to functions truncated by $\beta_\theta(D')$,
then $WW^*$ may be replaced
by the integral kernel
$$
K(r,x';s,y')=
\mu\int e^{i\langle\zeta,x'-z\rangle-i\langle\zeta_{s,r},y'-z_{s,r}\rangle}
g\bigl(\mu^{\frac 12}(x'-z)\bigr)\,g\bigl(\mu^{\frac 12}(y'-z_{s,r})\bigr)
\,\beta_\theta(\zeta)\,dz\,d\zeta\,,
$$
where we use the shorthand notation
\begin{equation}\label{zsrdef}
(z_{s,r},\zeta_{s,r})=\Theta_{s,r}(z,\zeta)\,.
\end{equation}

The factors $\beta_\theta(D')$ in \eqref{Fbound} can now be ignored
(since they are bounded in the
desired norms), and we are reduced to establishing mapping properties for $K$.
We observe that \eqref{Fbound}, and hence
Theorem \ref{maintheorem6},
can be reduced to establishing the following pair of bounds:
\begin{equation}\label{kest1}
\sup_{r,s\in[0,\eps]}\Bigl\|\int K(r,x';s,y')\,f(y')\,dy'\Bigr\|_{L^2(\R^2)}\le
\|f\|_{L^2(\R^2)}\,,
\end{equation}
and
\begin{equation}\label{kest2}
\Bigl\|\int K(r,x';s,y')\,f(y')\,dy\Bigr\|_{L^\infty_{x_2}L^2_{x_3}(\R^2)}
\lesssim
\mu\theta\,\bigr(\,1+\mu\theta^2\,|r-s|\,\bigr)^{-\frac 12}\,
\|f\|_{L^1_{y_2}L^2_{y_3}(\R^2)}\,.
\end{equation}
To see this, note that interpolation yields the bound
$$
\Bigl\|\int K(r,x';s,y')\,f(y')\,dy\Bigr\|_{L^q_{x_2}L^2_{x_3}(\R^2)}\lesssim
\bigl(\mu\theta\bigr)^{1-\frac 2q}\,
\bigr(\,1+\mu\theta^2\,|r-s|\,\bigr)^{\frac 1q-\frac 12}\,
\|f\|_{L^{q'}_{y_2}L^2_{y_3}(\R^2)}\,.
$$
If $q\ge 6$, then
$$
\bigl(\mu\theta\bigr)^{1-\frac 2q}
\bigr(\,1+\mu\theta^2\,|s|\,\bigr)^{\frac 1q-\frac 12}
\le
\mu^{1-\frac 4q}\theta^{1-\frac 6q}|s|^{-\frac 2q}
=\mu^{2\delta(q)}\theta^{1-\frac 6q}|s|^{-\frac 2q}\,,
$$
and by the Hardy-Littlewood-Sobolev inequality,
$\| |s|^{-\frac 2q}*f\|_{L^{q'}(\R)}\lesssim \|f\|_{L^q(R)}$,
we obtain the following bound equivalent to \eqref{Fbound}
$$
\Bigl\|\int K(r,x';s,y')\,F(s,y')\,ds\,dy'\Bigr\|_{L^q_rL^q_{x_2}L^2_{x_3}(S)}
\lesssim
\mu^{2\delta(q)}\theta^{1-\frac 6q}\,
\|F\|_{L^{q'}_sL^{q'}_{y_2}L^2_{y_3}(S)}\,.
$$

The bound \eqref{kest1} follows immediately from the $L^2$
boundedness of $T_\mu$, and the fact that $\Theta_{s,r}$ preserves
the measure $dz\,d\zeta\,,$ so it remains to establish
\eqref{kest2}. We start by estimating the derivatives of the
Hamiltonian flow with respect to the initial parameters $z$ and
$\zeta$. We only need bounds for curves lying entirely in the
region
$$\zeta_3\in[\tfrac 14\mu,2\mu]\quad
\text{and}
\quad \zeta_2 \in [\tfrac 14\theta\mu,2\theta\mu]
$$
(respectively $|\zeta_2|\le\mu^{\frac 12}$ in case
$\theta=\mu^{-\frac 12}$). In order to avoid extraneous powers of
$\mu$ it is convenient to exploit homogeneity to reduce to the
case $|\zeta|\approx 1$. For the purposes of the rest of this
section, we thus assume that the symbol $q$ (and hence the flow
$\Theta_{s,r}$) is homogeneous of degree one in $\zeta$, and agrees
with our previous definition of $q$ on the above region (which had
smoothly set $q=|\zeta|$ outside the region $|\zeta_3|\approx\mu\,,\,
|\zeta_2|\lesssim\frac 12\mu$.)

\begin{theorem}\label{flowtheorem}
Let $z_{s,r}$ and $\zeta_{s,r}$ be defined as functions of $(z,\zeta)$
by \eqref{zsrdef}.
Let $d_\zeta$ and $d_z$ respectively
denote the $\zeta$-gradient and $z$-gradient operators.
Then, for $\zeta_3=1$ and
$\zeta_2\approx\theta$ (respectively $|\zeta_2|\le\mu^{-\frac 12}$ in case
$\theta=\mu^{-\frac 12}$)
the following bounds hold.
\begin{align}
\bigl|d_z z_{s,r}-I\bigr|\lesssim |s-r|&\,,\qquad\qquad
\bigl|d_\zeta z_{s,r}\bigr|\lesssim |s-r|\,,
\notag\\
\label{estone}\\
\bigl|d_\zeta \zeta_{s,r}-I\bigr|\lesssim |s-r|&\,,\qquad\qquad
\bigl|d_z \zeta_{s,r}\bigr|\lesssim 1\,.\notag
\end{align}
Also,
\begin{align}
\bigl|d^2_z z_{s,r}\bigr|\lesssim
\langle \mu^{\frac 12}|s-r|&\,\rangle\,,\qquad\qquad
\bigl|d^2_z \zeta_{s,r}\bigr|\lesssim \mu^{\frac 12}\,
\notag\\
\label{estthree}\\
\bigl|d_z d_\zeta z_{s,r}\bigr|\lesssim
|s-r|\,\langle\mu^{\frac 12}|s-r|&\,\rangle\,,
\qquad\qquad
\bigl|d_z d_\zeta \zeta_{s,r}\bigr|\lesssim
\langle\mu^{\frac 12}|s-r|\,\rangle\,.
\notag
\end{align}
Furthermore, for $k\ge 2$,
\begin{equation}\label{esttwo}
\bigl|d^k_\zeta z_{s,r}|\;+\;
\bigl|d^k_\zeta \zeta_{s,r}|\;\lesssim
\; |s-r|\,\langle\mu^{\frac 12}|s-r|\,\rangle^{k-1}\,.
\end{equation}
\end{theorem}
\demo
We start with the relation
$$
z_{s,r}=z+\int_r^s (d_\zeta q)(t,z_{t,r},\zeta_{t,r})\,dt\,,\qquad
\zeta_{s,r}=\zeta-\int_r^s (d_z q)(t,z_{t,r},\zeta_{t,r})\,dt\,.
$$
Differentiating with respect to $z$ and $\zeta$ yields
\begin{equation}\label{intform1}
\begin{pmatrix}d z_{s,r}\\ \\
d\zeta_{s,r} \end{pmatrix}
=
\begin{pmatrix}dz\\ \\
d\zeta\end{pmatrix}+
\int_r^s M(t,z_{t,r},\zeta_{t,r})\cdot
\begin{pmatrix}d z_{t,r}\\ \\
d\zeta_{t,r} \end{pmatrix}\,dt\,,
\end{equation}
where
$$
M=\begin{pmatrix} \bigl(d_z d_\zeta q \bigr)
& \bigl( d_\zeta d_\zeta q\bigr)
\\ \\
-\bigl(d_z d_z q\bigr) & -\bigl(d_\zeta d_z q\bigr)
\end{pmatrix}
$$
The key estimate is that, for $i+j=2$,
$$
\int_r^s\bigl|(d^i_z d^j_\zeta q)(t,z_{t,r},\zeta_{t,r})\bigr|\,dt\lesssim
\begin{cases} |s-r|\,,&i\le 1\,,\\
\quad 1\,,&i=2\,.
\end{cases}
$$
This follows by \eqref{qsymbolcond} (recall
that $z_2$ equals $x_2$), and the fact
that $|(d_tz_{t,r})_2|\approx \theta$ in case $\theta>\mu^{-\frac 12}$.
In the case $\theta=\mu^{-\frac 12}$,
estimate \eqref{qsymbolcond} shows that the integrand is in
fact uniformly bounded.

An application of the Gronwall lemma yields
$$
|d z_{t,r}|\lesssim 1\,,\qquad |d\zeta_{t,r}|\lesssim 1\,,
$$
and plugging this into \eqref{intform1} yields \eqref{estone}.

To control higher order derivatives we proceed by induction.
For $k\ge 2$, we write
\begin{equation*}
\begin{pmatrix}d^k_\zeta z_{s,r}\\ \\
d^k_\zeta \zeta_{s,r} \end{pmatrix}
=
\int_r^s M(t,z_{t,r},\zeta_{t,r})\cdot
\begin{pmatrix}d^k_\zeta z_{t,r}\\ \\
d^k_\zeta \zeta_{t,r} \end{pmatrix}\,dt
+\int_r^s \begin{pmatrix}E_1(t) \\ \\ E_2(t)\end{pmatrix}\,dt
\end{equation*}
where $E_1(t)$ is a sum of terms of the form
$$
\bigl(d_z^{\,i} d_\zeta^{\,j+1} q\bigr)(t,z_{t,r},\zeta_{t,r})\cdot
\bigl(d^{k_1}_\zeta z_{t,r}\bigr)\cdots
\bigl(d^{k_i}_\zeta z_{t,r}\bigr)
\bigl(d^{k_{i+1}}_\zeta \zeta_{t,r}\bigr)\cdots
\bigl(d^{k_{i+j}}_\zeta \zeta_{t,r}\bigr)
$$
and $E_2$ is similarly a sum of such terms,
but with $d_z^{\,i+1}d_\zeta^{\,j}q$.
In both cases, $k_n<k$ for each $n$,
and $\;k_1+\cdots+k_{i+j}=k\,.$
By induction we may thus assume that the
estimates \eqref{estone} and \eqref{esttwo} hold
for all terms arising in $E_1$ and $E_2$.

The bound \eqref{qsymbolcond} implies, as above, that for $|\zeta|=1$
\begin{equation}\label{pintest}
\int_r^s |(d^{i+1}_z d^j_\zeta q)(t,z_{t,r},\zeta_{t,r})|\,dt\lesssim
\begin{cases}
|s-r|\,,&\quad i=0\\
\mu^{\frac 12(i-1)}\,,&\quad i\ge 1
\end{cases}
\end{equation}
The induction hypothesis yields that
$$
\bigl|
\bigl(d^{k_1}_\zeta z_{t,r}\bigr)\cdots
\bigl(d^{k_i}_\zeta z_{t,r}\bigr)
\bigl(d^{k_{i+1}}_\zeta \zeta_{t,r}\bigr)\cdots
\bigl(d^{k_{i+j}}_\zeta \zeta_{t,r}\bigr)
\bigr|
\lesssim |t-r|^i\,\langle\mu^{\frac 12}|t-r|\,\rangle^{k-i-j}\,.
$$
Together these yield
$$
\int_r^s |E_2(t)|\,dt\lesssim
|s-r|\,\langle\mu^{\frac 12}|s-r|\,\rangle^{k-1}\,,
$$
and the same holds for $E_1$.
The estimate \eqref{esttwo} follows by the Gronwall lemma.

To establish the first line of \eqref{estthree}, we write
\begin{equation}\label{intform2}
\begin{pmatrix}d^2_z z_{s,r}\\ \\
d^2_z \zeta_{s,r} \end{pmatrix}
=
\int_r^s M(t,z_{t,r},\zeta_{t,r})\cdot
\begin{pmatrix}d^2_z z_{t,r}\\ \\
d^2_z \zeta_{t,r} \end{pmatrix}\,dt
+\int_r^s \begin{pmatrix}E_1(t) \\ \\ E_2(t)\end{pmatrix}\,dt\,,
\end{equation}
where now
$$
\int_r^s |E_1(t)|\,dt\lesssim 1\,,\qquad
\int_r^s |E_2(t)|\,dt\lesssim \mu^{\frac 12}\,.
$$
A first application of Gronwall yields
$|d^2_z z_{s,r}|+|d^2_z \zeta_{s,r}|\lesssim \mu^{\frac 12}$,
and plugging this into \eqref{intform2} and
using \eqref{pintest} yields the first line of \eqref{estthree}.
The second line follows
by similar considerations.
\findemo

\begin{corollary}\label{dzetazest}
The following bounds hold for $\zeta_3=1$ and $\zeta_2\approx\theta$,
$$
\left|\,d_\zeta z_{s,r}- \int_r^s d^2_\zeta
q(t,\Theta_{t,r}(z,\zeta))\,dt\,\right| \le c\,|s-r|^2\,,
$$
where $c$ can be made small by taking the constant $c_0$ in
condition \eqref{cond3} small.
\end{corollary}
\demo
Given $c$, choosing the constant $c_0$ small yields the bounds
$$
|d_\zeta d_{z,\zeta}q|\le c\,,\qquad\qquad \int_r^s
|d_z^2q(t,\Theta_{t,r}(z,\zeta))|\,dt\le c\,.
$$
Together with the bounds \eqref{estone}, plugging this into
\eqref{intform1} yields
successively the bounds
$$
|d_\zeta\zeta_{s,r}-I|\le c\,|s-r|\,,\qquad\qquad \left|\,d_\zeta
z_{s,r}-\int_r^s d^2_\zeta q(t,\Theta_{t,r}(z,\zeta))\,dt \,\right|
\le c\,|s-r|^2\,.\qed
$$

\begin{lemma}
Suppose that $|\zeta|\approx \mu\,,$ and $\thetabar$ is a number with
$\thetabar\ge\mu^{-\frac 12}$ and $\mu\thetabar^2|s-r|\le 1\,.$ Then, for all
$\alpha$ and $j$,
\begin{equation}\label{symbest1}
\bigl|
\bigl(\zeta\cdot d_\zeta\bigr)^j
\bigl(\mu\thetabar\partial_\zeta\bigr)^\alpha
\mu^{\frac 32}\thetabar\, d_\zeta z_{s,r}\bigr|\lesssim 1\,,
\end{equation}
and for all $\alpha$ and $j$ with $j+|\alpha|\ge 1$,
\begin{equation}\label{symbest2}
\bigl|
\bigl(\zeta\cdot d_\zeta\bigr)^j\bigl(\mu\thetabar\partial_\zeta\bigr)^\alpha
\mu\thetabar \langle d_\zeta \zeta_{s,r},y-z_{s,r}\rangle
\bigr|\lesssim \langle\mu^{\frac 12}|y-z_{s,r}|\,\rangle\,.
\end{equation}
\end{lemma}
\demo
First consider \eqref{symbest1}.
By homogeneity of $z_{s,r}$ and its derivatives,
it suffices to consider $j=0$. We then have, by \eqref{estone}--\eqref{esttwo}
and homogeneity,
$$
\bigl|\bigl(\mu\thetabar\partial_\zeta\bigr)^\alpha
\mu^{\frac 32}\thetabar \,d_\zeta z_{s,r}\bigr|
\lesssim
\mu^{\frac 12}\thetabar^{|\alpha|+1}|s-r|\,
\langle\mu^{\frac 12}|s-r|\,\rangle^{|\alpha|}
\le \langle\mu^{\frac 12}\thetabar\,|s-r|\,\rangle^{|\alpha|+1}\lesssim 1\,.
$$
Here we use that $\mu\thetabar^2\ge 1$, so that
$\mu^{\frac 12}\thetabar|s-r|\le 1$.
For \eqref{symbest2}, note that if $|\alpha|=0$
and $j\ne 0$ then the term vanishes by homogeneity,
so we may assume $|\alpha|\ge 1$.
By homogeneity we may also restrict to the case $j=0$.
First consider the case where all derivatives fall on $\zeta_{s,r}$.
The resulting term is bounded by
$$
\mu\thetabar^{|\alpha|+1}|s-r|\,
\langle\mu^{\frac 12}|s-r|\,\rangle^{|\alpha|}\,|y-z_{s,r}|
\lesssim
\langle\mu^{\frac 12}\thetabar\,|s-r|\,\rangle^{|\alpha|}
\langle\mu^{\frac 12}|y-z_{s,r}|\,\rangle
\lesssim \langle\mu^{\frac 12}|y-z_{s,r}|\,\rangle\,.
$$
If one or more derivatives falls on $z_{s,r}$, the term is bounded by
$$
\mu\thetabar^{|\alpha|+1}|s-r|\,
\langle\mu^{\frac 12}|s-r|\,\rangle^{|\alpha|-1}
\lesssim \mu\thetabar^2|s-r|\,
\langle\mu^{\frac 12}\thetabar\,|s-r|\,\rangle^{|\alpha|-1}
\lesssim 1\,.
$$
\findemo

Recall that the kernel we are proving \eqref{kest2} for is
$$
K(r,x';s,y')=
\mu\int e^{i\langle\zeta,x'-z\rangle-i\langle\zeta_{s,r},y'-z_{s,r}\rangle}
g\bigl(\mu^{\frac 12}(x'-z)\bigr)\,g\bigl(\mu^{\frac 12}(y'-z_{s,r})\bigr)
\,\beta_\theta(\zeta)\,dz\,d\zeta\,,
$$
where $\beta_\theta(\zeta)$ is a cutoff to $\zeta_3\approx\mu$ and
$\zeta_2\approx\theta\mu$, respectively
$|\zeta_2|\lesssim\mu^{\frac 12}$ in case
$\theta=\mu^{-\frac 12}$.

In what follows, for the case $\theta>\mu^{-\frac 12}$
we will need to consider finer angular
decompositions in $\zeta$, depending on $|s-r|$.
We will assume, for the following theorem,
that $\beta_\thetabar(\zeta)$ is a smooth cutoff to
a set of the form
$$
\zeta_3\approx\mu\,,
\qquad
\zeta_2\approx\theta\mu\,,
\qquad
\Bigl|\,\frac{\zeta_2}{\zeta_3}-\theta'\,\Bigr|\lesssim\thetabar\,,
$$
where $\theta'\approx\theta$, and where $\mu^{-\frac 12}
\le\thetabar\le\theta$.
For $\theta=\mu^{-\frac 12}$, we need consider only
$\beta_\thetabar=\beta_\theta$.

\begin{theorem}\label{ktheorem}
Consider the kernel $K$ with $\beta_\theta(\zeta)$ replaced by
$\beta_\thetabar(\zeta)$, with $\beta_\thetabar$ as above. Suppose
that $\mu\thetabar^2|s-r|\le 1$. Fix a vector $\xi'$ in the
support of $\beta_\thetabar(\zeta)$, and let
$(x'_{s,r},\nu_{s,r})$ be the projection of $\Theta_{s,r}(x',\xi')$
onto the cosphere bundle. Thus, $x'_{s,r}=z_{s,r}$ and
$\nu_{s,r}=|\zeta_{s,r}|^{-1}\zeta_{s,r}$ if $z=x'$ and
$\zeta=\xi'$. Then
$$
|K(r,x';s,y')|\lesssim
\mu^2\thetabar\,
\bigl(\,1+\mu\thetabar\,|\,y'-x'_{s,r}|+
\mu\,|\,\langle\nu_{s,r},y'-x'_{s,r}\rangle\,|\,\bigr)^{-N}
\,.
$$
\end{theorem}
\demo
We introduce the differential operators, where $z_{s,r}$ and $\zeta_{s,r}$
are as in \eqref{zsrdef},
$$
L_1=
\frac{1-i\bigl(\,\langle\zeta,x'-z\rangle
-\langle\zeta_{s,r},y'-z_{s,r}\rangle
\,\bigr)
\,\langle\zeta,d_\zeta\rangle}
{1+\bigl|\,\langle\zeta,x'-z\rangle-\langle\zeta_{s,r},y'-z_{s,r}\rangle
\,\bigr|^2}
\,,
$$
and
$$
L_2=
\frac{1-i\mu\thetabar\,\bigl(\,x'-z-d_\zeta\zeta_{s,r}
\cdot(y'-z_{s,r})\,\bigr)\cdot d_\zeta}
{1+\mu^2\thetabar^2\,\bigl|\,x'-z-d_\zeta\zeta_{s,r}\cdot(y'-z_{s,r})\bigr|^2}
\,.
$$
Each of these preserves the phase function in $K$,
and an integration by parts argument,
using the estimates \eqref{symbest1} and \eqref{symbest2},
bounds $|K(r,x';s,y')|$ by
the following integral
\begin{multline*}
\mu\int
\bigl(\,1+\mu\thetabar\,\bigl|\,x'-z-d_\zeta\zeta_{s,r}\cdot(y'-z_{s,r})
\bigr|\,\bigr)^{-N}
\bigl(\,1+\bigl|\,\langle\zeta,x'-z\rangle-
\langle\zeta_{s,r},y'-z_{s,r}\rangle\,\bigr|\,
\bigr)^{-N}\\
\times
\bigl(\,1+\mu^{\frac 12}|x'-z|\,\bigr)^{-N}
\bigl(\,1+\mu^{\frac 12}|y'-z_{s,r}|\,\bigr)^{-N}\,dz\,d\zeta\,,
\end{multline*}
where the integral is over the support of $\beta_\thetabar(\zeta)$,
which has volume
$\mu^2\thetabar\,.$ 
We will show below that
\begin{multline}\label{zetaineq}
\mu\thetabar\,|\,d_\zeta\zeta_{s,r}\cdot (x'_{s,r}-z_{s,r})-(x'-z)|+
|\langle\zeta_{s,r},x'_{s,r}-z_{s,r}\rangle-\langle\zeta,x'-z\rangle|\\
\lesssim 1+\mu\,|x'-z|^2\,.
\end{multline}
This implies that the integrand is dominated by
$$
\bigl(\,1+\mu\thetabar\,|d_\zeta \zeta_{s,r}\cdot(y'-x'_{s,r})|+
|\langle\zeta_{s,r},y'-x'_{s,r}\rangle|
\,\bigr)^{-N}
\bigl(\,1+\mu^{\frac 12}|x'-z|\,\bigr)^{-N}\,.
$$
By \eqref{estone}, the matrix $d_\zeta \zeta_{s,r}$ is invertible. Also
by \eqref{estone}, the angle of $\zeta_{s,r}$ to $\mu\nu_{s,r}$ is less than
$\thetabar+|x'-z|$. Since $\mu\thetabar\ge\mu^{\frac 12}$,
and $|\zeta_{s,r}|\approx\mu$, together these
dominate the integrand by
$$
\bigl(\,1+\mu\thetabar\,|\,y'-x'_{s,r}|+
\mu\,|\langle\nu_{s,r},y'-x'_{s,r}\rangle|
\,\bigr)^{-N}
\bigl(\,1+\mu^{\frac 12}|x'-z|\,\bigr)^{-N}\,,
$$
from which the theorem follows easily.

We now establish \eqref{zetaineq}.
Consider the first term on the left. By homogeneity, we may assume that
$1=|\zeta|=|\xi'|$, so that $|\zeta-\xi'|\le\thetabar\,.$
By \eqref{estthree} and Taylor's theorem, we then have
\begin{multline*}
|x'_{s,r}-z_{s,r}-(d_z z_{s,r})(x'-z)-(d_\zeta z_{s,r})(\xi'-\zeta)|
\\
\lesssim \langle \mu^{\frac 12}|s-r|\,\rangle\,|x'-z|^2+
|s-r|\,\langle\mu^{\frac 12}|s-r|\,
\rangle\bigl(\thetabar^2+\thetabar|x'-z|\,\bigr)\,.
\end{multline*}
After multiplication by $\mu\thetabar$, each term on the right is
bounded by $1+\mu\,|x'-z|^2\,.$
Also,
$$
\mu\thetabar\,|\,(d_\zeta z_{s,r})(\zeta-\xi')|\lesssim
\mu\thetabar^2\,|s-r|\le 1\,.
$$
Since $d\zeta_{s,r}\wedge dz_{s,r}=d\zeta\wedge dz$, we have
$$
\partial_{\zeta_i} \zeta_{s,r}\cdot \partial_{z_j} z_{s,r}
-\partial_{\zeta_i} z_{s,r}\cdot \partial_{z_j}\zeta_{s,r}=\delta_{ij}\,,
$$
where $\cdot$ pairs the $z_{s,r}$ and $\zeta_{s,r}$ indices.
By \eqref{estone}, we have
$$
\mu\thetabar\,\bigl|d_\zeta z_{s,r}\bigr|\,\bigl|d_z\zeta_{s,r}\bigr|\,|x'-z|
\lesssim
\mu^{\frac 12}|x'-z|\,.
$$
Together, this yields
$$
\mu\thetabar\,|\,d_\zeta\zeta_{s,r}\cdot(x'_{s,r}-z_{s,r})-(x'-z)|
\lesssim 1+\mu\,|x'-z|^2\,,
$$
which concludes the bound for the first term.

To handle the second term, it suffices by homogeneity to show that,
for $|\zeta|=|\xi'|=1$,
$$
|\langle\zeta_{s,r},x'_{s,r}-z_{s,r}\rangle-\langle\zeta,x'-z\rangle|
\lesssim |x'-z|^2+\thetabar^2|s-r|\,.
$$
We calculate
\begin{multline*}
\frac{d}{ds}\langle\zeta_{s,r},x'_{s,r}-z_{s,r}\rangle=\\
-\bigl\langle (d_z
q)(s,\Theta_{s,r}(z,\zeta)),x'_{s,r}-z_{s,r}\bigr\rangle+
\bigl\langle\zeta_{s,r},(d_\zeta q)(s,\Theta_{s,r}(x',\xi'))-
(d_\zeta q)(s,\Theta_{s,r}(z,\zeta))\bigr\rangle\,.
\end{multline*}
By homogeneity, the right hand side equals
\begin{equation}\label{perror}
q(s,\Theta_{s,r}(x',\xi'))-q(s,\Theta_{s,r}(z,\zeta))-
\bigl(\Theta_{s,r}(x',\zeta)-\Theta_{s,r}(z,\zeta)\bigr)\cdot
(d_{z,\zeta}q)(s,\Theta_{s,r}(z,\zeta))
\end{equation}
plus an error which, since $q_\zeta$ is Lipschitz, is bounded by
\begin{equation}\label{chibound}
|\Theta_{s,r}(x',\xi')-\Theta_{s,r}(z,\zeta)|^2\lesssim
|x'-z|^2+\thetabar^2\,.
\end{equation}
Let
$\gamma_\sigma(t)=\sigma\Theta_{s,r}(x',\xi')
+(1-\sigma)\Theta_{s,r}(z,\zeta)\,.$
Then \eqref{perror} equals
$$
\int_0^1\int_r^s(1-\sigma)\bigl(\Theta_{s,r}(x',\xi')
-\Theta_{s,r}(z,\zeta)\bigr)^2
(d_{z,\zeta}^2 q)(t,\gamma_\sigma(t))\,dt\,d\sigma\,.
$$
By \eqref{chibound}, the integral of terms involving
$d_zd_\zeta q$ and $d_\zeta^2 q$
are bounded by
$$
|s-r|\,|x'-z|^2+|s-r|\,\thetabar^2\le|x'-z|^2+\thetabar^2|s-r|\,.
$$
The integral of terms involving $d_z^2 q$ are bounded by
$$
\Bigl(\;\sup_{r\le t\le s}|x'_{t,r}-z_{t,r}|^2\,\Bigr)\;
\sup_\sigma\int_r^s|(d^2_z q)(t,\gamma_\sigma(t))|\,dt
\lesssim |x'-z|^2+\thetabar^2|s-r|^2\,,
$$
where we use \eqref{estone}, \eqref{qsymbolcond}, and the fact that
$(\dot{\gamma}_\sigma)_2\approx\theta$ in the case 
$\theta>\mu^{-\frac 12}$.
\findemo

\noindent{\bf Proof of estimate \eqref{kest2}.}
We establish \eqref{kest2} by showing that
\begin{equation}\label{kest3}
\sup_{x_2,x_3,y_2}\int\,|K(r,x';s,y')|\,dy_3\lesssim
\mu\theta\,\bigl(\,1+\mu\theta^2|s-r|\,\bigr)^{-\frac 12}\,.
\end{equation}
Transposing $(s,y')$ and $(r,x')$ in the formula for $K$
leads to the same kernel
if $\beta_\theta(\zeta)$ is replaced by 
$\beta_\theta(\zeta_{r,s}(y',\zeta))$,
and the same proof will show that
$$
\sup_{x_2,y_2,y_3}\int\,|K(r,x';s,y')|\,dx_3\lesssim
\mu\theta\,\bigl(\,1+\mu\theta^2|s-r|\,\bigr)^{-\frac 12}\,,
$$
yielding \eqref{kest2} by Schur's lemma.

Suppose first that $\mu\theta^2|s-r|\le 1$.
Then \eqref{kest3} follows immediately
from Theorem \ref{ktheorem} with $\thetabar=\theta$,
since $\nu_{s,r}=|\zeta_{s,r}|^{-1}\zeta_{s,r}$
is within a small angle of the $\xi_1$ axis.

If $\mu\theta^2|s-r|>1$, we let
$\thetabar=\mu^{-\frac 12}|s-r|^{-\frac 12}$, and
decompose $K$ into a sum of terms by writing
$\beta_\theta(\zeta)=\sum_j\beta_j(\zeta)$,
with each $\beta_j(\zeta)$ a cutoff to a sector of angle $\thetabar$.

We fix $\eta^j$ in the support of $\beta_j(\zeta)$,
with
$$
(\eta^j)_3=\mu\,,\quad\text{and}\quad
|(\eta^i)_2-(\eta^j)_2|\approx\mu\,\thetabar\,|i-j|\,.
$$
We then have decomposed
$K=\sum_j K_j$, where by Theorem \eqref{ktheorem}
$$
|K_j(r,x';s,y')|\lesssim
\mu^2\thetabar\,
\bigl(\,1+
\mu\thetabar\,|\,y'-w^j_{s,r}|+\mu\,|\langle\nu^j_{s,r},y'-w^j_{s,r}\rangle|
\,\bigr)^{-N}
\,,
$$
where $(w^j_{s,r},\nu^j_{s,r})$ is the projection onto the
cosphere bundle of $\Theta_{s,r}(x',\eta^j)$. Since
$(\nu^j_{s,r})_3\approx 1$, we have
$$
\int |K_j(r,x';s,y')|\,dy_3\lesssim
\mu\thetabar\,\bigl(1+\mu\thetabar\,|\,y_2-(w^j_{s,r})_2|\,\bigr)^{-N}\,.
$$
Since $\mu\thetabar\approx\mu\theta\,(1+\mu\theta^2|s-r|\,)^{-\frac 12}$,
it suffices
to show that
$$
\sup_{x_2,x_3,y_2}
\sum_j\bigl(1+\mu\thetabar\,|\,y_2-(w^j_{s,r})_2|\,\bigr)^{-N}\lesssim 1\,,
$$
which we do by recalling that $\mu\thetabar^2|s-r|=1$, and
showing that
$$
|(w^i_{s,r})_2-(w^j_{s,r})_2|\approx \thetabar\,|s-r|\,|i-j|\,.
$$
We finally show this by noting that, for $\zeta_3=1$ and
$|\zeta_2|\le \tfrac 12$,
we have $d_{\zeta_2}^2 q\approx 1\,.$
Corollary \ref{dzetazest} thus yields
$d_{\zeta_2} (z_{s,r})_2 \approx s-r$ for such $\zeta$.
Consequently,
$$
|z_{s,r}(x',\eta^i)-z_{s,r}(x',\eta^j)|\approx
\mu^{-1}\,|s-r|\,|\,(\eta^i)_2-(\eta^j)_2|\approx
\thetabar\,|s-r|\,|i-j|\,.\qed
$$

\newsection{Energy Flux Estimates}\label{section6}
In this section we complete the proof of Theorem
\ref{maintheorem3} by establishing the endpoint estimates where
$q=8$.
 We do this by establishing the nested square-summability condition
\eqref{cijcond}. Recall that we are assuming
$$
D_1 u_\l-P_\l u_\l=F_\l\,,
$$
where $2P_\l=p_\l(x,D')+p_\l(x,D')^*$, and we write
$$
D_1 u_j-P_j u_j=F_j+G_j\,,
$$
where $u_j=\beta_j(D')u_\l$, the operator $P_j=p_j(x,D')$ has
symbol truncated to $x'$- frequencies
less than $\l^{\frac 12}\theta_j^{-\frac 12}$, and
\begin{align}\label{6.1}
F_j&=\beta_j(D')F_\l+[\beta_j(D'),P_j]u_\l+
\beta_j(D')\bigl(P_\l-p_\l(x,D')\bigr)u_\l\,,
\\
\label{6.2}
G_j&=\beta_j(D')\bigl(p_\l(x,D')-p_j(x,D')\bigr)u_\l\,.
\end{align}
Let
\begin{multline}\label{cjkdef}
c_{j,k}=\|u_j\|_{L^\infty L^2(S_{j,k})}+
\l^{\frac 14}\theta_j^{\frac 14}
\|\langle \l^{\frac 12}\theta_j^{-\frac 12}x_2\rangle^{-1}u_j\|_{L^2(S_{j,k})}
\\
+\|F_j\|_{L^1L^2(S_{j,k})}
+\l^{-\frac 14}\theta_j^{-\frac 14}
\|\langle\l^{\frac 12}\theta_j^{-\frac 12}x_2\rangle^2 G_j\|_{L^2(S_{j,k})}\,.
\end{multline}
We need to show that
\begin{equation}\label{cijcond'}
\sum_{j=1}^{N_\l}
c_{j,k(j)}^2\lesssim \|u_\l\|_{L^\infty L^2(S)}^2+\|F_\l\|_{L^2(S)}^2\,,
\end{equation}
where $k(j)$ denotes any sequence
of values for $k$ such that the slabs $S_{j,k(j)}$ are nested,
in that for $j\ge 1$ we have
$S_{j+1,k(j+1)}\subset S_{j,k(j)}$.
The analogous bound for $j<0$ will follow by an identical proof.

\subsection{Estimates on $u_j$}
We begin by establishing the square-summability
estimates for the first two terms on the
right hand side of \eqref{cjkdef}.
By translation invariance we may assume each
$S_{j,k(j)}$ contains $x_1=0$. We then take
$S_j$ to be the slab $[0,\eps\,2^{-j}]\times\R^2$,
and will show that
$$
\sum_i
\Bigl(\;\|u_i\|^2_{L^\infty L^2(S_i)}
+\l^{\frac 12}\theta_i^{\frac 12}
\|\langle\l^{\frac 12}\theta_i^{-\frac 12}x_2\rangle^{-1}u_i\|^2_{L^2(S_i)}
\Bigr)
\lesssim
\|u_\l\|_{L^\infty L^2}^2+\|F_\l\|_{L^2}^2\,.
$$
The same bounds will hold for $x_1\in[-\eps\,2^{-j},0]$.

Since $F_\l\in L^1_{x_1}L^2_{x'}$, by Duhamel we can reduce matters
to the homogeneous case $F_\l=0$. Assume this, and let $f(x')=u_\l(0,x')$.
Let $W$ denote the solution operator for the Cauchy problem
associated to $P_\l$, so that $u_\l=W\! f$,
It then suffices to show that
\begin{multline}\label{Wfijest}
\|\beta_i(D')W\beta_j(D')f\|_{L^\infty L^2(S_i)}
+
\l^\frac 14\theta_i^{\frac 14}\|
\langle\l^{\frac 12}\theta_i^{-\frac 12}x_2\rangle^{-1}
\beta_i(D')W\beta_j(D')f\|_{L^2(S_i)}
\\
\lesssim 2^{-\frac 34|i-j|}\,\|f\|_{L^2}\,.
\end{multline}
To prove \eqref{Wfijest}, we will construct for each given $j$
a function $v$ which satisfies the following conditions.
\begin{equation}\label{vcond1}
v(0,x')=\beta_j(D')f(x')\,,
\qquad \beta_i(D)v=0\quad\text{if}\quad |i-j|\ge 5\,,
\end{equation}
\begin{align}\label{vcond2}
\|v\|_{L^\infty L^2(S_j)}&\lesssim \|f\|_{L^2}
\\
\label{vcond3}
\l^{\frac 14}\theta_j^{\frac 14}\|
\langle\l^{\frac 12}\theta_j^{-\frac 12}x_2\rangle^{-1}
v\|_{L^2(S_j)}& \lesssim \|f\|_{L^2}\,,
\end{align}
and such that
\begin{align}\label{vcond4}
\|D_1v-P_\l v\|_{L^1 L^2(S_j)}&\lesssim
(\l^{\frac 12}\theta_j^{\frac 32})^{-\frac 12}\,\|f\|_{L^2}\,,
\\
\label{vcond5}
\|D_1v-P_\l v\|_{L^1 L^2(S_j)}&\lesssim
(\l^{\frac 12}\theta_j^{\frac 32})^{-1}\,
\|\langle\l^{\frac 12}\theta_j^{-\frac 12}x_2\rangle f\|_{L^2}\,.
\end{align}

Let us show that these imply the estimate \eqref{Wfijest}.
Consider the first term on the left hand side of \eqref{Wfijest}.
We will prove the stronger statement
\begin{equation}\label{LinfWest}
\|\beta_i(D')W(r)\beta_j(D')f\|_{L^2_{x'}}
\lesssim 2^{-\frac 34|i-j|}\,\|f\|_{L^2}\,, \quad |r|\le
\eps\,\max(2^{-i},2^{-j})\,.
\end{equation}
By self adjointness
(the adjoint of $W(r)$ is the wave map going the other way),
we can then assume that $\theta_j=2^{-j}\ge \theta_i=2^{-i}$.
This assumption
now means we need to control data at angle $2^{-j}$ for time $\eps\,2^{-j}$.

We write $W\beta_j(D')f=v-w$. The desired estimate holds for
the $v$ term by \eqref{vcond2}, since we may assume $|i-j|\le 4$
by \eqref{vcond1} (and we may shrink $\eps$ by a factor of 16.)

To control $w$, we note that
$$
w(0,x')=0\,,\qquad D_1w-P_\l w=D_1v-P_\l v\,.
$$
Energy estimates and \eqref{vcond4} thus yield
\begin{equation}\label{Linfwest}
\|w\|_{L^\infty L^2(S_j)}\le
(\l^{\frac 12}\theta_j^{\frac 32})^{-\frac 12}\,\|f\|_{L^2}\,.
\end{equation}
Since $\l^{\frac 12}\ge 2^{\frac 32 i}$, this yields the desired
bound on $u$.

To estimate the second term in \eqref{Wfijest}, we first consider
the case $\theta_i\le \theta_j$.
We again write $W\beta_j(D')f=v-w$, and note that the desired estimate on
$v$ follows by \eqref{vcond3} and \eqref{vcond1}.
(The operator $\beta_i(D')$ preserves the $L^2$-weight
$\langle\l^{\frac 12}\theta_i^{-\frac 12}x_2\rangle^{-1}$
since $\lambda^{\frac 12}\theta_i^{-\frac 12}\le 2^{-i}\l$.)

The estimate on $w$ for $\theta_i\le \theta_j$ follows by \eqref{Linfwest},
$$
\l^\frac 14\theta_i^{\frac 14}\|w\|_{L^2(S_i)}
\le
\l^\frac 14\theta_i^{\frac 34}\|w\|_{L^\infty L^2(S_j)}
\lesssim
\theta_j^{-\frac 34}\theta_i^{\frac 34}\|f\|_{L^2}\,.
$$

Now consider the case $\theta_j\le\theta_i$.
The above steps handle the case $|i-j|\le 4$, so we assume $i\ge j+5$.
We take adjoints to reduce matters to showing that, for $j\ge i+5$,
\begin{equation*}
\Bigl\|\;\int_{|s|\le \eps \theta_i}\beta_j(D')W(s)^*\beta_i(D')
F(s,\cd)\,ds\,\Bigr\|_{L^2}
\\
\lesssim
\l^{-\frac 14}\theta_i^{-\frac 14}\,2^{-\frac 34|i-j|}
\|\langle\l^{\frac 12}\theta_i^{-\frac 12} x_2\rangle F\|_{L^2(S_i)}\,.
\end{equation*}
This bound, in turn,
follows from showing that, for $|r|\le \eps\, 2^{-i}$ and $j\ge i+5$,
$$
\|\beta_j(D')W^*(r)\beta_i(D')
f\|_{L^2}
\lesssim
(\l^{\frac 12}\theta_i^{\frac 32})^{-1}
\|\langle\l^{\frac 12}\theta_i^{-\frac 12} x_2\rangle f\|_{L^2}\,.
$$
We may replace $W^*(r)$ by $W(r)$, since $W^*$ is the Cauchy map for
data at $x_1=r$ to $x_1=0$, and after exchanging $i$ and $j$ this bound
is a consequence of \eqref{vcond5}.

\subsection{The construction of $v$.}
We assume that $\theta_j$ is now fixed, and rescale spatial variables
by $\theta_j$. We thus need to construct $v$ on the slab
$S=[0,\eps]\times\R^2$.
As before, let $\mu=\l\theta_j$, and let
$\beta_j(D')$ denote the rescaled localization operators, which
will localize to $\xi_2\approx\theta_j\mu\,,$ $\,\xi_3\approx\mu$.
Let $f$ denote the rescaled initial data $\beta_j(D')f(\theta_j\cd)$.

In these rescaled variables
it suffices to produce $v$ satisfying
\begin{equation}\label{vcond1'}
v(0,x')=f(x')\,,
\qquad \beta_i(D)v=0\quad\text{if}\quad |i-j|\ge 5\,,
\end{equation}
\begin{align}\label{vcond2'}
\|v\|_{L^\infty L^2(S)}&\lesssim \|f\|_{L^2}
\\
\label{vcond3'}
\mu^{\frac 14}\theta_j^{\frac 12}\|
\langle\mu^{\frac 12}x_2\rangle^{-1}
v\|_{L^2(S)}& \lesssim \|f\|_{L^2}\,,
\end{align}
and such that
\begin{align}\label{vcond4'}
\|D_1v-Q_\mu v\|_{L^2(S)}&\lesssim
(\mu^{\frac 12}\theta_j)^{-\frac 12}\,\|f\|_{L^2}\,,
\\
\label{vcond5'}
\|D_1v-Q_\mu v\|_{L^1L^2(S)}&\lesssim
(\mu^{\frac 12}\theta_j)^{-1}\,
\|\langle\mu^{\frac 12}x_2\rangle f\|_{L^2}\,.
\end{align}
Here $Q_\mu$ is the rescaled operator $P_\l$, which has symbol
truncated to $x'$-frequencies less than $c\mu$.

We will construct $v$ using the modified FBI/C\'ordoba-Fefferman transform
$T_\mu$ introduced in \S\ref{section4}. The key idea is that this
transform conjugates the operator $Q_\mu$ to the Hamiltonian
flow field, plus a bounded error which is roughly local. Precisely,
we will show that
$$
T_\mu Q_\mu T_\mu^* = D_q+K,
$$
where $D_q$ is the Hamiltonian vector field of the symbol $q$ (which
we recall is
frequency localized to $\mu^{\frac 12}$), and where $K$ is an operator
on $L^2_{x',\xi'}$, depending on parameter $x_1$, for which
we establish weighted $L^2$ estimates.

The transform $\tilde u=T_\mu u$
of the exact solution $u$ to $D_1u-Q_\mu u=0$, with initial
data $f$, satisfies
$$
D_1\tilde u-D_q\tilde u
=K\tilde u\,,\qquad \tilde u(0,x',\xi')=\tilde f(x',\xi')\,.
$$
The operator $K$ will introduce terms which are well-behaved
after integration along the flow of $D_1-D_q$ at angle $\theta_j$.
We will construct the approximate solution
$v$ by truncating the operator $K$ to such angles. For this purpose
we introduce cutoffs $\phi_j(\xi')$ and $\psi_j(\xi')$,
with slightly larger supports than $\beta_j(\xi')$,
such that
\begin{align*}
\text{dist}\bigl(\text{supp}(1-\phi_j),\text{supp}(\beta_j)\bigr)&
\ge 2^{-j-10}\mu\,,\\
\text{dist}\bigl(\text{supp}(1-\psi_j),\text{supp}(\phi_j)\bigr)&
\ge 2^{-j-10}\mu\,,
\end{align*}
and also that
$$
\text{dist}\bigl(\text{supp}(\psi_j),\text{supp}(\beta_i)\bigr)
\ge 2^{-j-10}\mu \quad\text{if}\quad|i-j|\ge 5\,.
$$

The $\xi'$-support of $\tilde f$ lies in the $c\mu^{-\frac 12}$
neighborhood of the support of $\beta_j(\xi')$.
Since $c\ll 1$,
$\theta_j\ge \mu^{-\frac 12}$, and
$|d_x q(x,\xi')|\le c\,\theta_j\,|\xi'|\,,$
we can assume that every integral curve of $D_1-D_q$ passing through
this neighborhood
remains $\xi'$-distance at least $2^{-10}\mu\theta_j$
away from the support of $(1-\phi_j)$.

Furthermore, we can assume that
any integral curve of $D_1-D_q$ passing at any point
through the support of $\psi_j$
does not meet the $c\mu^{-\frac 12}$ neighborhood
of the support of $\beta_i(\xi')$, provided $|i-j|\ge 5$.

We will take $v=T_\mu^* \tilde v$ where $\tilde v$ solves
\begin{equation}\label{Lveqn}
D_1\tilde v-D_q\tilde v=\psi_j\,K\tilde v\,,\qquad
\tilde v(0,x',\xi')=\tilde f(x',\xi')\,.
\end{equation}
The cutoff $\psi_j$ restricts the right hand side
to $\xi_2\approx\theta_j\mu$,
where the integral of $K$ along $D_1-D_q$ is under control.
Furthermore, since the support of $\tilde v$ will be contained
in the union of the integral curves of $D_1-D_q$ passing through the support
of $\psi_j$ at some point $x_1$, then $v$ will
satisfy $\beta_i(D')v=0$ for $|i-j|\ge 5$.

Next, since $Q_\mu T_\mu^*=T_\mu^*D_q+T_\mu^*K$, it holds that
$$
D_1v-Q_\mu v=-T_\mu^*\bigl((1-\psi_j)K\tilde v\bigr)\,,
$$
so estimates \eqref{vcond4'} and \eqref{vcond5'} will follow from
\begin{equation}\label{vcond4''}
\|(1-\psi_j)\,K\tilde v\|_{L^2(\tilde S)}
\lesssim
(\mu^{\frac 12}\theta_j)^{-\frac 12}\,\|\tilde f\|_{L^2}\,,
\end{equation}
and
\begin{equation}\label{vcond5''}
\|(1-\psi_j)\,K\tilde v\|_{L^1L^2(\tilde S)}
\lesssim
(\mu^{\frac 12}\theta_j)^{-1}
\|\langle\mu^{\frac 12}x_2\rangle \tilde f\|_{L^2}\,.
\end{equation}
where $\tilde S=[0,\eps]\times\R^4_{x',\xi'}\,.$

We thus need to show that $K\tilde v$ is small away from the
set $\xi_2\approx \theta_j\mu$, which we do by establishing
weighted norm estimates on $\tilde v$, and
decay estimates on the kernel $K$.
The weights involve the natural distance function on $\R^4_{x',\xi'}$
associated to the C\'ordoba-Fefferman transform,
$$
\text{dist}_\mu(x',\xi';y',\eta')
=\mu^{\frac 12}\,|x'-y'|+\mu^{-\frac 12}|\xi'-\eta'|\,.
$$
Let $K(x',\xi';y',\eta')$ denote the integral kernel of $K$ (we supress
the parameter $x_1$). Then we will show that
\begin{multline}\label{kest}
|K(x',\xi';y',\eta')|\lesssim
\bigl(\,1+{\text{dist}}_\mu(x',\xi';y',\eta')\bigr)^{-N}\\
+
c_0\,\mu^{\frac 12}\theta_j\langle\, \mu^{\frac 12}x_2\rangle^{-N}
\langle\,\mu^{-\frac 12}|\xi_2-\eta_2|\,\rangle^{-2}\,
\bigl(\,1+\mu^{\frac 12}|x'-y'|+\mu^{-\frac 12}|\xi_3-\eta_3|\,\bigr)^{-N}\,,
\end{multline}
where $c_0$ is the small constant of \eqref{cond3}.

Let $E_0$ be the subset of $R^4_{x',\xi'}$
$$
E_0=\R^2_{x'}\times \text{supp}(\beta_j(\xi'))\,,
$$
and let $E_{x_1}$ be the image of $E_0$ under the flow along $D_1-D_q$
for time $x_1$. We consider the weight function
$$
M(x,\xi')=M_{x_1}(x',\xi')=1+\text{dist}_\mu(x',\xi';E_{x_1})\,.
$$

The weighted norm estimates we establish for solutions of \eqref{Lveqn} are
\begin{equation}\label{Mest1}
\|M\tilde v\|_{L^\infty L^2(\tilde S)}\lesssim \|M\!\tilde f\|_{L^2}\,,
\end{equation}
\begin{equation}\label{Mest2}
\|\langle\mu^{\frac 12}x_2\rangle^{-1} M \tilde v\|_{L^2(\tilde S)}
\lesssim (\mu^{\frac 12}\theta_j)^{-\frac 12}\|M\!\tilde f\|_{L^2}\,,
\end{equation}
and
\begin{equation}\label{Mest3}
\|\langle\mu^{\frac 12}x_2\rangle^{-2} M \tilde v\|_{L^1L^2(\tilde S)}
\lesssim (\mu^{\frac 12}\theta_j)^{-1}
\|\langle\mu^{\frac 12}x_2\rangle M\!\tilde f\|_{L^2}\,.
\end{equation}

Let us show how \eqref{vcond2'}--\eqref{vcond5'}
follow from \eqref{kest} and
\eqref{Mest1}--\eqref{Mest3}.
The bounds \eqref{vcond2'} and \eqref{vcond3'} are direct
consequences of \eqref{Mest1} and \eqref{Mest2}, since $M=1$ on the support
of $\tilde f$. Also,
\eqref{vcond4'}--\eqref{vcond5'} follow from
\eqref{vcond4''}--\eqref{vcond5''}, so we focus on
\eqref{vcond4''}--\eqref{vcond5''}.

We write $K=K_1+K_2$, where the kernels $K_1$ and $K_2$ are respectively
dominated by the first and second terms on the right hand side of
\eqref{kest}.

First note that, since
$\text{dist}_{\xi'}(\supp(1-\phi_j),E_{x_1})
\ge 2^{-10}\mu\theta_j$ for all $x_1$,
it follows from \eqref{Mest1}--\eqref{Mest3} that
\begin{align*}
\|(1-\phi_j)\tilde v\|_{L^\infty L^2(\tilde S)}
&\lesssim
(\mu^{\frac 12}\theta_j)^{-1}\|\tilde f\|_{L^2}\,,\\
\|\langle\mu^{\frac 12}x_2\rangle^{-1}(1-\phi_j)
\tilde v\|_{L^2(\tilde S)}
&\lesssim
(\mu^{\frac 12}\theta_j)^{-\frac 32}\|\tilde f\|_{L^2}\,,\\
\|\langle\mu^{\frac 12}x_2\rangle^{-2}(1-\phi_j)
\tilde v\|_{L^1L^2(\tilde S)}
&\lesssim
(\mu^{\frac 12}\theta_j)^{-2}
\|\langle\mu^{\frac 12}x_2\rangle\tilde f\|_{L^2}\,.
\end{align*}
By the bounds on $K_1$ and $K_2$ and Schur's Lemma we thus have
\begin{align*}
\|K_1(1-\phi_j)\tilde v\|_{L^\infty L^2(\tilde S)}&\lesssim
(\mu^{\frac 12}\theta_j)^{-1}\|\tilde f\|_{L^2}\\
\|K_2(1-\phi_j)\tilde v\|_{L^2(\tilde S)}&\lesssim
(\mu^{\frac 12}\theta_j)^{-\frac 12}\|\tilde f\|_{L^2}\,,\\
\|K_2(1-\phi_j)\tilde v\|_{L^1L^2(\tilde S)}&\lesssim
(\mu^{\frac 12}\theta_j)^{-1}
\|\langle \mu^{\frac 12}x_2\rangle \tilde f\|_{L^2}\,.
\end{align*}
Next, we note that the integral of $K_1$, as well as the integral
of $(\mu^{\frac 12}\theta_j)^{-1}\langle\mu^{\frac 12}y_2\rangle K_2$,
over the set $|\xi'-\eta'|\ge 2^{-10}\mu\theta_j$ is bounded
by $(\mu^{\frac 12}\theta_j)^{-1}$, which yields by
\eqref{Mest1}--\eqref{Mest3} that
\begin{align*}
\|(1-\psi_j)K_1\phi_j\tilde v\|_{L^\infty L^2(\tilde S)}&\lesssim
(\mu^{\frac 12}\theta_j)^{-1}\|\tilde f\|_{L^2}\\
\|(1-\psi_j)K_2\phi_j\tilde v\|_{L^2(\tilde S)}&\lesssim
(\mu^{\frac 12}\theta_j)^{-\frac 12}\|\tilde f\|_{L^2}\,,\\
\|(1-\psi_j)K_2\phi_j\tilde v\|_{L^1L^2(\tilde S)}&\lesssim
(\mu^{\frac 12}\theta_j)^{-1}
\|\langle\mu^{\frac 12}x_2\rangle\tilde f\|_{L^2}\,.
\end{align*}
Together these yield the estimates \eqref{vcond4''} and \eqref{vcond5''}.

We turn to the proof of estimates \eqref{Mest1}--\eqref{Mest3}.
\begin{lemma}\label{KMlemma}
Take $E\subset\R^4$ and let $M(x',\xi')= 1+\dist_\mu(x',\xi';E)$.
Also, let $r_-=\frac 12(|r|-r).$
Then, for postive integers $k$ and $n$, and real number $r$,
$$
\|M \langle\mu^{\frac 12}x_2\rangle^k
\langle\mu^{\frac 12}(x_2-r)_-\rangle^n
 K_1 g\|_{L^2(\R^4)}
\lesssim
\|M \langle \mu^{\frac 12}x_2\rangle^k
\langle\mu^{\frac 12}(x_2-r)_-\rangle^n
g\|_{L^2(\R^4)}\,,
$$
and
\begin{multline*}
\|M \langle \mu^{\frac 12}x_2\rangle^k
\langle\mu^{\frac 12}(x_2-r)_-\rangle^n K_2 g\|_{L^2(\R^4)}\\
\lesssim
c\,\mu^{\frac 12}\theta_j
\|M \langle \mu^{\frac 12}x_2\rangle^{k-N}
\langle\mu^{\frac 12}(x_2-r)_-\rangle^n g\|_{L^2(\R^4)}\,.
\end{multline*}
The bounds are uniform over all subsets $E\subset \R^4$ and real numbers $r$.
\end{lemma}
\demo
Let $K_0$ denote the integral kernel
$$
K_0(x',\xi';y',\eta')=\bigl(\,1+\mu^{-\frac 12}|\eta_2-\xi_2|\,\bigr)^{-2}
\bigl(\,1+\mu^{\frac 12}|y'-x'|+\mu^{-\frac 12}|\eta_3-\xi_3|\,\bigr)^{-N}\,.
$$
By the rapid decrease of $K$ in $x'$ and $\xi_3$,
both estimates are a simple of the following bound
$$
\|M K_0g\|_{L^2}\lesssim \|M g\|_{L^2}\,.
$$
By making the measure preserving change of variables
$(x',\xi')\rightarrow(\mu^{\frac 12}x',\mu^{-\frac 12} \xi')$,
we may assume $\mu=1$.
By the rapid decrease of $K_0$ in the $x'$ and $\xi_3$ variables, we may
bound
\begin{multline*}
\|\,M(x',\xi')\int K_0(x',\xi';y',\eta')\,g(y',\eta')\,dy'\,d\eta'\,
\|_{L^2(dx'\,d\xi')}
\\
\lesssim
\|\,\int M(y',\xi_2,\eta_3) \langle\xi_2-\eta_2\rangle^{-2}\,
g(y',\eta')\,d\eta_2\,\|_{L^2(d\xi_2\,dy'\,d\eta_3)}\,.
\end{multline*}
We lastly use the following consequence of the Calder\'on
commutator theorem \cite{Ca}.
\begin{lemma}
Let $M(r)$ denote a weight function on the real line, satisfying
$$
M(r)\ge 1\,,\qquad |M(r)-M(s)|\le |r-s|\,.
$$
Then the convolution kernel $\langle r\rangle^{-2}$ is bounded
on $L^2(M(r)dr)$ by a uniform constant.
\end{lemma}
\demo
We need to show that the integral kernel
$$
\frac{M(r)\,M(s)^{-1}}{\langle r-s\rangle^2}=
\frac{M(r)-M(s)}{\langle r-s\rangle^2}\,M(s)^{-1}
+\frac{1}{\langle r-s\rangle^2}
$$
is bounded on $L^2(dr)$. Since $M(s)^{-1}\le 1$
and the latter kernel
is integrable, it suffices to show that the map
$$f\;\rightarrow\;\int_{-\infty}^\infty \frac{M(r)-M(s)}{\langle \, r-s\,
\rangle^2} \, f(s)\, ds$$ is bounded on $L^2(dr)$.  Clearly
$$f\;\rightarrow\;\int_{|r-s|\le 1} \frac{M(r)-M(s)}{\langle \, r-s\,
\rangle^2} \, f(s)\, ds$$ is bounded on $L^2(dr)$, and so it
suffices to show that
\begin{equation}\label{Cald}
\left\| \, \int_{|r-s|> 1} \frac{M(r)-M(s)}{\langle \, r-s\,
\rangle^2} \, f(s)\, ds\, \right\|_{L^2(dr)}\le C\|f\|_{L^2(dr)}.
\end{equation}
But
$$\left| \, \frac{M(r)-M(s)}{\langle \, r-s\,
\rangle^2}-\frac{M(r)-M(s)}{(r-s)^2}\, \right|\le
\frac{|M(r)-M(s)|}{\langle \, r-s\, \rangle^2(r-s)^2}\le
\frac1{|r-s|^3},$$ which means that \eqref{Cald} holds if and only
if the map
$$f\;\rightarrow\;\int_{|r-s|> 1} \frac{M(r)-M(s)}{(
r-s)^2} \, f(s)\, ds$$ is bounded on $L^2(dr)$.  But since $M$ is
Lipschitz, this follows from the classical commutator
estimate of Calder\'on (Theorem 2 in \cite{Ca}).
%the kernel
%$$
%\frac{M(r)-M(s)}{\langle r-s\rangle^2}
%$$
%is bounded on $L^2(dr)$. This is the Calderon commutator theorem, which
%is usually stated
%with $\langle r-s\rangle^{-2}$ replaced by $(r-s)^{-2}$. Indeed,
%in our case we may assume $|M''(r)|\lesssim 1$, since we can smooth
%out $M$ at scale 1, in which case the difference between the two
%kernels is the Hilbert transform plus an integrable kernel.
\findemo

In the following steps, we will use $r$ and $s$ as real variables
that take the place of $x_1$.

Let $J_{r,s}:\R^4\rightarrow \R^4$ denote the flow along $D_1-D_q$, starting
at the slice $x_1=s$ and ending at $x_1=r$. We will also use
$J_{r,s}$ to denote the unitary map on $L^2(\R^4)$
$$
\bigl(J_{r,s}f\bigr)(x',\xi')=f\bigl(J_{s,r}(x',\xi')\bigr)\,.
$$
This map is unitary since the Hamiltonian flow is symplectic, hence
preserves $dx'\,d\xi'$.
We also use the fact that,
if $|\xi_2|\,,\,|\eta_2|\approx\mu\theta_j$
and $|\xi_3|\,,\,|\eta_3|\approx \mu$, then
the map $J_{r,s}$ approximately preserves $\dist_\mu$, in that
$$
\dist_\mu\bigl(J_{r,s}(x',\xi');J_{r,s}(y',\eta')\bigr)
\approx
\dist_\mu(x',\xi';y',\eta')\,.
$$
By homogeneity of the Hamiltonian flow, this follows from the fact that
the flow is Lipschitz on the set $|\xi'|=1$, which is a consequence
of Theorem \ref{flowtheorem}.

The function $\tilde v$ satisfies
$$
D_1\tilde v-D_q\tilde v
=\psi_j K\tilde v\,,\qquad \tilde v(0,x',\xi')=\tilde f(x',\xi')
\,.
$$
Let $\U$ denote the map, taking the space of
functions on $\tilde S$ to itself,
defined by
$$
\U F(r,\cd)=\int_0^rJ_{r,s}\psi_j KF(s,\cd)\,ds\,.
$$
Thus, $(D_1-D_q)\,\U F=\psi_j K F\,.$ If we let $F(r,\cd)=J_{r,0}\tilde f\,,$
so that $D_1F-D_q F=0$, then we can formally write the solution $\tilde v$ as
$$
\tilde v=\sum_{n=0}^\infty \U^n F\,.
$$
We need to show this sum converges in the appropriate norm, which
we do by showing that $\U$ is a contraction. We split $\,\U=\U_1+\U_2$,
corresponding to the splitting $K=K_1+K_2$.

The estimates we require for $\,\U_1$ are:
\begin{equation}\label{Uest1}
\|M \U_1 F\|_{L^\infty L^2(\tilde S)}
\lesssim
\eps\,\|M F\|_{L^\infty L^2(\tilde S)}
\end{equation}
\begin{equation}\label{Uest2}
(\mu^{\frac 12}\theta_j)^{\frac 12}
\|M\langle\mu^{\frac 12}x_2\rangle^{-1}
\U_1 F\|_{L^2(\tilde S)}\lesssim
\eps\,\|M F\|_{L^\infty L^2(\tilde S)}\,.
\end{equation}

For the $\,\U_2$ term we require the bounds:
\begin{equation}\label{Uest4}
\|M\U_2 F\|_{L^\infty L^2(\tilde S)}
\lesssim
c\,(\mu^{\frac 12}\theta_j)^{\frac 12}
\|M \langle\mu^{\frac 12}x_2\rangle^{-1}
F\|_{L^2(\tilde S)}
\end{equation}
\begin{equation}\label{Uest5}
\|M\langle\mu^{\frac 12}x_2\rangle^{-1}\U_2
F\|_{L^2(\tilde S)}
\lesssim
c\,\|M \langle\mu^{\frac 12}x_2\rangle^{-1}
F\|_{L^2(\tilde S)}
\end{equation}

The inequality \eqref{Uest1}
is a consequence of Lemma \ref{KMlemma} with $k=n=N=0$, and the
fact that $J_{r,s}$ preserves the distance function $\dist_\mu$,
hence the weight $M$.

For \eqref{Uest2}, we apply Cauchy-Schwarz to yield
$$
|M \U_1 F|^2(r,x',\xi')\lesssim
\eps\,\int_0^\eps \bigl| (M\psi_jK_1F)(s,J_{s,r}(x',\xi'))\bigr|^2\,ds\,.
$$
We multiply by $\langle\mu^{\frac 12}x_2\rangle^{-2}$ and integrate
$dx'\,d\xi'$, changing variables by $J_{s,r}$ on the right, to obtain
\begin{multline*}
\|M\langle\mu^{\frac 12}x_2\rangle^{-1}
\U_1 F\|^2_{L^2(\tilde S)}\\
\lesssim
\eps\,\int_0^\eps\int_0^\eps\int_{\R^4}
\langle\mu^{\frac 12}x_2\circ J_{r,s}\rangle^{-2}
\bigl|M\psi_jK_1F\bigr|^2(s,x',\xi')\,dx'\,d\xi'\,ds\,dr\,.
\end{multline*}
We next observe that, for $\xi'$ in the support of $\psi_j$,
\begin{equation}\label{gammaintest}
\int\langle\mu^{\frac 12}x_2\circ J_{r,s}\rangle^{-2}\,dr\lesssim
(\mu^{\frac 12}\theta_j)^{-1}\,,
\end{equation}
which holds since $\frac{dx_2}{dr}\approx\theta_j$.
Lemma \ref{KMlemma} with $k=n=N=0$
now yields \eqref{Uest2}.

To show \eqref{Uest4}, we write
\begin{align*}
|(M\U_2 F)(r,x',\xi')|^2&\lesssim
\left|\int_0^r
\bigl(M \psi_j K_2 F\bigr)(s,J_{s,r}(x',\xi'))\,ds\,\right|^2\\
&\lesssim
(\mu^{\frac 12}\theta_j)^{-1}
\int_0^\eps
\bigl|M \langle\mu^{\frac 12}x_2\rangle
\psi_j K_2 F\bigr|^2 (s,J_{s,r}(x',\xi'))\,ds
\end{align*}
where we use \eqref{gammaintest}.
To conclude \eqref{Uest2} we take the integral $dx'\,d\xi'$ of both
sides, using the fact that $J_{s,r}$ preserves the measure, and
apply Lemma \ref{KMlemma} with $k=1$, $n=0$, and $N=2$.

For \eqref{Uest5}, we write as above
$$
|(M\U_2 F)(r,x',\xi')|^2
\lesssim
(\mu^{\frac 12}\theta_j)^{-1}
\int_0^\eps
\bigl|M \langle\mu^{\frac 12}x_2\rangle
\psi_j K_2 F\bigr|^2 (s,J_{s,r}(x',\xi'))\,ds\,.
$$
For $\xi'$ in the support of $\psi_j$ we have
$$
\langle\mu^{\frac 12}x_2\rangle^{-2}\,
\langle\mu^{\frac 12}x_2\circ J_{s,r}\rangle^{-2}\lesssim
\langle\mu^{\frac 12}\theta_j|r-s|\,\rangle^{-2}\,,
$$
and consequently
\begin{multline*}
\langle\,\mu^{\frac 12}x_2\rangle^{-2}\,
|(M\U_2 F)(r,x',\xi')|^2\\
\lesssim
(\mu^{\frac 12}\theta_j)^{-1}
\int_0^\eps
\langle\,\mu^{\frac 12}\theta_j|r-s|\,\rangle^{-2}
\bigl|M \langle\mu^{\frac 12}x_2\rangle^2
\psi_j K_2 F\bigr|^2 (s,J_{s,r}(x',\xi'))\,ds\,.
\end{multline*}
We take the integral $dx'\,d\xi'$, changing variables by $J_{s,r}$ on the
right, and apply Lemma \ref{KMlemma}
with $k=2$, $n=0$ and $N=3$, to yield \eqref{Uest5}.

We now turn to the proof of \eqref{Mest1}--\eqref{Mest3}.
First, note that by \eqref{Uest1}--\eqref{Uest2} and
\eqref{Uest4}--\eqref{Uest5}, for small $c$ and $\eps$ the map $\U$
is a contraction in the norm
$$
|||F|||=
\|MF\|_{L^\infty L^2(\tilde S)}
+\mu^{\frac 14}\theta_j^{\frac 12}
\|M\langle\mu^{\frac 12}x_2\rangle^{-1} F\|_{L^2(\tilde S)}\,.
$$
Recall that
$\tilde v=\sum_{n=0}^\infty \U^n F\,,$ where
$F(r,\cd)=T_{r,0}\tilde f\,.$
Furthermore,the bound \eqref{gammaintest} yields
$$|||F|||\lesssim \|M\tilde f\|_{L^2}\,.$$
Consequently
$$
|||\tilde v|||\lesssim \|M\tilde f\|_{L^2}\,,
$$
which implies \eqref{Mest1} and \eqref{Mest2}.

To derive \eqref{Mest3}, we use the fact that each of the estimates
\eqref{Uest1}--\eqref{Uest5} holds if $M$ is in each instance
replaced by the weight
$$
M\langle\mu^{\frac 12}(x_2-c_2\theta_jx_1)_-\rangle\,,
$$
where $c_2>0$ is a constant such that
$\frac{dx_2}{dr}> c_2\theta_j$ on
curves of $D_1-D_q$ passing through the support of $\psi_j$.
This holds since
$\langle\mu^{\frac 12}(\cd)_-\rangle$ is positive and decreasing,
and hence, if $\psi_j(\xi')\ne 0$ and $r\ge s$, then
$$
\langle\mu^{\frac 12}(x_2\circ J_{r,s}-c_2\theta_j r)_-\rangle
\le
\langle\mu^{\frac 12}(x_2-c_2\theta_j s)_-\rangle\,.
$$

Consequently,
\begin{align*}
\|M\langle\mu^{\frac 12}(x_2-c_2\theta_jx_1)_-\rangle
\langle\mu^{\frac 12}x_2\rangle^{-1}\tilde v\|_{L^2(\tilde S)}
&\lesssim
(\mu^{\frac 12}\theta_j)^{-\frac 12}
\|M\langle\mu^{\frac 12}(x_2)_-\rangle\tilde f\|_{L^2}\\
&\lesssim
(\mu^{\frac 12}\theta_j)^{-\frac 12}
\|M\langle\mu^{\frac 12}x_2\rangle\tilde f\|_{L^2}\,.
\end{align*}
On the other hand, for $x_1>0$,
$$
\langle\mu^{\frac 12}(x_2-c_2\theta_jx_1)_-\rangle\,
\langle\mu^{\frac 12}x_2\rangle\gtrsim
\langle\mu^{\frac 12}\theta_j x_1\rangle\,,
$$
hence
$$
\|M\langle\mu^{\frac 12}x_2\rangle^{-2}\tilde v\|_{L^1 L^2(\tilde S)}
\lesssim
(\mu^{\frac 12}\theta_j)^{-\frac 12}
\|M\langle\mu^{\frac 12}(x_2-c_2\theta_jx_1)_-\rangle
\langle\mu^{\frac 12}x_2\rangle^{-1}\tilde v\|_{L^2(\tilde S)}
\,,
$$
yielding \eqref{Mest3}.

\subsection{The estimate on $K$}
We establish here the estimate \eqref{kest} for the integral kernel $K$
defined by
$$
T_\mu Q_\mu T_\mu^*=D_q+K\,.
$$
Here, $Q_\mu=\frac 12\bigl(q_\mu(x,D')+q_\mu(x,D')^*\bigr)$, where
$q_\mu$ is the symbol $p_\l$ rescaled by $\theta_j$, and hence
truncated to $x'$-frequencies less than $c\mu$. The symbol $q$, on the
other hand, is
obtained by truncating $q_\mu$ to $x'$-frequencies less than $c\mu^{\frac 12}$.

It is a simple consequence of Lemma \ref{pest} and Lemma \ref{wavetype}
that the kernel of the operator
$$
T_\mu\,q(x,D')^*\,T_\mu^*-D_q
$$
satisfies the estimate \eqref{kest}.
By taking adjoints the same applies with $q(x,D')^*$ replaced by $q(x,D')$,
and we are reduced to establishing the estimates \eqref{kest} for the
kernel of the operator
$$
T_\mu\,(q_\mu(x,D')-q(x,D')\,)\,T_\mu^*\,.
$$

The kernel
$K(x',\xi';y',\eta')$ of this operator takes the form
(we suppress the irrelevant parameter $x_1$)
\begin{equation*}
\int
e^{i\langle \zeta',z'-y'\rangle}\,
e^{-i\langle \xi',z'-x'\rangle}\,
\Bigl[q_\mu(z',\zeta')-q(z',\zeta')\Bigr]\,
\widehat{g}\bigl(\mu^{-\frac 12}(\zeta'-\eta')\bigr)\,
g\bigl(\mu^{\frac 12}(z'-x')\bigr)\,dz'\,d\zeta'\,.
\end{equation*}

Suppose that $p(x')$ is a smooth function on $x_2\ge 0$, which is
constant for $x_2\ge 1$. We extend $p$ in an even manner to $x_2\le 0$.
Let $q_\mu=S_\mu[p(\theta_j\cd)]$, and
$q=S_{\!\sqrt{\mu}}[p(\theta_j\cd)]\,,$
where $S_\l$ denotes smooth truncation of the Fourier transform to
frequencies less than $c\l$.
It then follows that
\begin{align}\label{qmucond}
\bigl|D_{x'}^\beta\bigl(q_\mu(x')-q(x')\bigr)\bigr|
&\lesssim
\theta_j\mu^{\frac 12(|\beta|-1)}\langle\mu^{\frac 12}x_2\rangle^{-N}\,
\|D_{x'}p\|_{C^N(x_2\ge 0)}\,,\quad
|\beta|\le 1\,,\\
\notag
\bigl|D_{x'}^2\bigl(q_\mu(x')-q(x')\bigr)\bigr|
&\lesssim
\theta_j\Bigl(\mu^{\frac 12}\langle\mu^{\frac 12}x_2\rangle^{-N}+
\mu\langle\mu x_2\rangle^{-N}\Bigr)
\|D_{x'}p\|_{C^{N+2}(x_2\ge 0)}\,.
\end{align}
Indeed, it suffices to verify these bounds for $x_2>0$, and by splitting
up $p$ to separately consider the case that $p$ is smooth across $x_2=0$
and constant for $|x_2|\ge 1$, and the case that $p$ is smooth on $x_2\le 0$
and vanishes for $x_2\ge 0$. The latter case is handle by simple size bounds
on the convolution kernels $S_\mu$ and $S_{\!\sqrt{\mu}}$.
For the smooth part, we have bounds
\begin{align*}
|p_\l-S_{\theta_j^{-1}\mu^{\frac 12}}p_\l|(\theta_j x')
&\lesssim
(\theta_j^{-1}\mu^{\frac 12})^{-1-N}\langle\theta_j x_2\rangle^{-N}
\,\|D_{x'}^{N+1}p\|_{C^0}\\
&\lesssim
\theta_j\,\mu^{-\frac 12}\langle\mu^{\frac 12}x_2\rangle^{-N}
\,\|D_{x'}^{N+1}p\|_{C^0}\,,
\end{align*}
and the same bounds apply to derivatives.

By the condition \eqref{cond3}, we easily obtain the following bounds
for $|\zeta'|\approx\mu$ and $\beta_2\le 2$,
$$
\bigl|\partial_{z'}^\beta\partial_{\zeta'}^{\alpha}\bigl(\,
q_\mu(z',\zeta')-q(z',\zeta')\bigr)\bigr|
\lesssim
c_0\,\theta_j\Bigl(\mu^{\frac 12}\langle\mu^{\frac 12}z_2\rangle^{-N}+
\mu\langle\mu z_2\rangle^{-N}\Bigr)
\mu^{\frac 12 |\beta_2|-|\alpha|}\,.
$$
In the formula for $K$ we can integrate by parts at will with respect to
$\mu^{\frac 12}D_{\zeta'}$ and $\mu^{-\frac 12}D_{z_3}$, and
twice with respect to $\mu^{-\frac 12}D_{z_2}$, to dominate $K$ by
\begin{multline*}
c_0\,\theta_j\int
\Bigl(\mu^{\frac 12}\langle\mu^{\frac 12}z_2\rangle^{-N}+
\mu\langle\mu z_2\rangle^{-N}\Bigr)
\langle\mu^{\frac 12}(z'-x')\rangle^{-N}
\langle\mu^{\frac 12}(z'-y')\rangle^{-N}\\
\times\langle\mu^{-\frac 12}(\zeta'-\eta')\rangle^{-N}
\langle\mu^{-\frac 12}(\zeta_3-\xi_3)\rangle^{-N}
\langle\mu^{-\frac 12}(\zeta_2-\xi_2)\rangle^{-2}
\,dz'\,d\zeta'
\end{multline*}
which is dominated by
$$
c_0\,\theta_j\,\mu^{\frac 12}\langle\mu^{\frac 12}x_2\rangle^{-N}\,
\langle\mu^{\frac 12}(x'-y')\rangle^{-N}\,
\langle\mu^{-\frac 12}(\xi_3-\eta_3)\rangle^{-N}\,
\langle\mu^{-\frac 12}(\xi_2-\eta_2)\rangle^{-2}
$$
yielding the desired bounds on $K$.\qed

We note here that similar considerations to the above yield the bounds,
for $|\xi'|\approx\mu$,
\begin{equation}\label{qsymbolcond'}
\bigl|\partial_x^\beta\partial^\alpha_{\xi'} q(x,\xi')\bigr|
\lesssim\begin{cases}
\mu^{1-|\alpha|}\,,\qquad |\beta|=0\,,\\
c_0\,
\bigl(\,1+\mu^{\frac 12(|\beta|-1)}\theta_j\,\langle
\mu^{\frac 12}x_2\rangle^{-N}\,\bigr)
\,\mu^{1-|\alpha|}\,,\quad|\beta|\ge 1\,.
\end{cases}
\end{equation}

\subsection{Estimates on $F_j$ and $G_j$.}
We conclude by establishing the square summability of the inhomogeneities
$F_j$ and $G_j$. Recall that
\begin{align*}
F_j&=\beta_j(D')F_\l+[\beta_j(D'),P_j]u_\l+
\beta_j(D')\bigl(P_\l-p_\l(x,D')\bigr)u_\l\,,
\\
G_j&=\beta_j(D')\bigl(p_\l(x,D')-p_j(x,D')\bigr)u_\l\,.
\end{align*}
We need to show that
$$
\sum_j\|F_j\|_{L^1L^2(S_j)}^2+
\l^{-\frac 14}\theta_j^{-\frac 14}
\|\langle\l^{\frac 12}\theta_j^{-\frac 12}x_2\rangle^2 G_j\|_{L^2(S_j)}^2
\lesssim
\|u_\l\|_{L^\infty L^2(S)}^2+\|F_\l\|_{L^2(S)}^2\,.
$$

The first term in $F_j$ is handled by noting that
$$
\sum_j\|\,\beta_j(D')F_\l\|_{L^1L^2(S_j)}^2
\le
\sum_j\|\,\beta_j(D')F_\l\|_{L^2(S)}^2\le\|F_\l\|_{L^2(S)}^2\,.
$$

Consider next the term $[P_j,\beta_j(D')]u_\l$. Since the symbol of
$P_j$ is truncated to $x'$-frequencies less than
$c\mu^{\frac 12}\le c\,\theta_j\mu$,
it holds that
$$
[P_j,\beta_j(D')]u_\l=[P_j,\beta_j(D')]\,\phi_j(D')u_\l\,.
$$
We claim that, uniformly over $x_1$,
\begin{equation}\label{Pjcommest}
\bigl\|[P_j,\beta_j(D')]f\bigr\|_{L^2_{x'}}\lesssim 2^j\,\|f\|_{L^2_{x'}}\,.
\end{equation}
Given this, we can bound
$$
\bigl\|[P_j,\beta_j(D')]\phi_j(D')u_\l\bigr\|_{L^1L^2(S_j)}\lesssim
2^j\,\|\,\phi_j(D') u_\l\|_{L^1L^2(S_j)}\le
\|\,\phi_j(D') u_\l\|_{L^\infty L^2(S_j)}\,,
$$
since $S_j$ is of length $2^{-j}$ in $x_1$. Since $\phi_j(D') u_\l$
involves $\beta_j(D')u_\l$ only for $|i-j|\le 4$, this term is square summable
over $j$.

To prove \eqref{Pjcommest}, it suffices to replace $P_j$ by $p_j(x,D')$.
The symbol $p_j$ equals $|\xi'|$ outside of the region $|\xi'|\approx\l$,
and $p_j(x',\xi')$ satisfies $S^1_{1,0}$ estimates for $x'$ derivatives
of order at most 1. Consquently, after subtracting off the term $|D'|$,
we may take $p_j(x,D')$ to have kernel $K_1(x',x'-y')$ where
$$
|K_1(x',z')|+|D_x'K_1(x',z')|\lesssim
\l\cdot\l^2\,\bigl(\,1+\l|z'|\,\bigr)^{-N}\,.
$$
On the other hand, $\beta_j(D')$ is a convolution kernel $K_2(x'-y')$
where
$$
\|z'\,K_2(z')\|_{L^1_{z'}}\lesssim 2^j\l^{-1}\,.
$$
The estimate \eqref{Pjcommest} follows by applying Taylor's theorem to
$$
\bigl[K_1,K_2\bigr](x',y')=
\int K_2(x'-z')\bigl(K_1(x',z'-y')-K_1(z',z'-y')\bigr)\,dz'\,.
$$
To control the last term in $F_j$ we note that, since the estimates on $K_1$
above also apply to $p_\l(x,D')$, we have uniform bounds
$$
\bigl\|\bigl(P_\l-p_\l(x,D')\bigr)f\bigr\|_{L^2_{x'}}
=\tfrac 12\bigl\|\bigl(p_\l(x,D')^*-p_\l(x,D')\bigr)f\bigr\|_{L^2_{x'}}
\lesssim\|f\|_{L^2_{x'}}\,.
$$
The last term in $F_j$ is orthogonal over $j$, and thus has square sum
bounded by $\|u\|_{L^2(S)}$.

We now estimate the term $G_j$. We split this up
$$
G_j=
\beta_j(D')\bigl(p_\l(x,D')-p_j(x,D')\bigr)\phi_j(D') u_\l+
\beta_j(D')p_\l(x,D')\bigl(1-\phi_j(D')\bigr)u_\l\,.
$$
Consider the second term in $G_j$. We write
$$
\beta_j(D')p_\l(x,D')\bigl(1-\phi_j(D')\bigr)u_\l=
\beta_j(D')\bigl(p_\l(x,D')-p_{\l\theta_j}(x,D')\bigr)
\bigl(1-\phi_j(D')\bigr)u_\l
$$
where $p_{\l\theta_j}$ is the symbol $p$ truncated to $x'$-frequencies
of size less than $c\l\theta_j$. The symbol
$p_\l-p_{\l\theta_j}$ is supported in the region $|\xi'|\approx\l$,
and by arguments similar to those deriving \eqref{qmucond}
(without the rescaling step), we have the estimates
$$
|\partial^\alpha_{\xi'}(p_\l-p_{\l\theta_j})(x,\xi')|\lesssim
\theta_j^{-1}\langle \l\theta_j x_2\rangle^{-N}\,\l^{-|\alpha|}\,.
$$
Its integral kernel is thus bounded by
$$
\theta_j^{-1}\langle \l\theta_j x_2\rangle^{-N}\,
\l^2\,\bigl(\,1+\l\,|x'-y'|\,\bigr)^{-N}\,.
$$
Since $\l\theta_j\ge \l^{\frac 12}\theta_j^{-\frac 12}$, it follows that,
uniformly in $x_1$,
$$
\l^{-\frac 14}\theta_j^{-\frac 14}
\|\langle\l^{\frac 12}\theta_j^{-\frac 12}x_2\rangle^2
\bigl(p_\l(x,D')-p_{\l\theta_j}(x,D')\bigr)f\|_{L^2_{x'}}\le
\l^{-\frac 14}\theta_j^{-\frac 54}\|f\|_{L^2_{x'}}\,,
$$
and the same holds for $p_\l(x,D')-p_{\l\theta_j}(x,D')$ replaced by
$\beta_j(D')\bigl(p_\l(x,D')-p_{\l\theta_j}(x,D')\bigr)$ since
$\beta_j(D')$ averages on scale smaller than
$\l^{-\frac 12}\theta_j^{\frac 12}$.
Thus
$$
\l^{-\frac 14}\theta_j^{-\frac 14}
\|\langle\l^{\frac 12}\theta_j^{-\frac 12}x_2\rangle^2\beta_j(D')
p_\l(x,D')\bigr(1-\phi_j(D')\bigr)u_\l\|_{L^2(S_j)}\le
\l^{-\frac 14}2^{\frac 34 j}\|u_\l\|_{L^\infty L^2(S)}\,.
$$
Since $2^j$ runs from 1 to $\l^{\frac 13}$, the right hand side
is square summable over $j$.

Recalling that the symbol $p_j(x,\xi')$ is truncated to $x'$-frequencies less
than $\l^{\frac 12}\theta_j^{-\frac 12}$, similar arguments
show that
\begin{multline*}
\l^{-\frac 14}\theta_j^{-\frac 14}
\|\langle\l^{\frac 12}\theta_j^{-\frac 12}x_2\rangle^2
\beta_j(D')\bigl(p_\l(x,D')-p_j(x,D')\bigr)\phi_j(D') u_\l
\|_{L^2(S_j)}\\
\lesssim
\l^{\frac 14}\theta_j^{\frac 14}
\|\langle\l^{\frac 12}\theta_j^{-\frac 12}x_2\rangle^{-1}\phi_j(D')u_\l
\|_{L^2(S_j)}\,.
\end{multline*}
The right hand side involves $u_i$ for $|i-j|\le 4$, hence, by the
earlier estimate for $u_i$, is square summable over $j$.

\newsection{Results for higher dimensions}\label{section7}

We show here that the steps of the preceeding sections
yield sharp $L^q$ estimates for spectral clusters on compact Riemannian
manifolds $M$ with boundary, of dimension $n\ge 3$, provided $q$
is sufficiently large. Precisely, we have the following

\begin{theorem}\label{penultimatemaintheorem}
Suppose that $u$ solves the Cauchy problem on $\R\times M$
\begin{equation*}
\partial_t^2u(t,x)=Pu(t,x)\,,
\qquad
u(0,x)=f(x)\, ,
\qquad \partial_tu(0,x)=0\,,
\end{equation*}
and satisfies either Dirichlet conditions
\begin{equation*}
 u(t,x)=0\quad\text{if}\quad x\in\partial M\,,
\end{equation*}
or Neumann conditions, where $N_x$ is a unit normal field with respect to $\g$,
\begin{equation*}
N_x\cdot\nabla_x u(t,x)=0\quad\text{if}\quad x\in\partial M\,.
\end{equation*}
Then the following bounds hold for $4\le q\le \infty$, if $n\ge 4$,
and $5\le q\le\infty$ if $n=3$.
\begin{equation*}
\|u\|_{L^q_xL^2_t(M\times[-1,1])}
\le C\,\|f\|_{H^{\delta(q)}(M)}\,,
\qquad \textstyle
\delta(q)=n\bigl(\frac 12-\frac 1q\bigr)-\frac 12\,.
\end{equation*}
\end{theorem}

These bounds of course imply that the estimates in \eqref{1.9}
hold for the spectral projector operators, $\chi_\lambda$, when
$q\ge 5$ for $n=3$ and $q\ge4$ if $n\ge 4$.

As noted in the introduction, these estimates are expected to hold 
in the larger range
$q\ge \frac{6n+4}{3n-4}$, in which case they (and their interpolation 
with the trivial
$L^2$ estimate) would be best possible. Establishing this larger 
range would require
exploiting dispersion in directions tangent to $\partial M$ for time 1, rather
than times on the order of the microlocalization angle $\theta$.

Following the earlier sections, we work in a neighborhood of $\partial M$ in
geodesic normal coordinates, and extend the operator $P$ evenly, 
and solution $u$
oddly or evenly, in the case of Dirichlet or Neumann conditions respectively.
We set $x_{n+1}=t$, and $x'=(x_2,\ldots,x_{n+1})$.

We then fix a frequency scale $\l$ and microlocalization angle 
$\theta_j\in[\l^{-\frac 13},c]\,.$
After factorizing $DA_\l D$, we set
$$
q(x,\xi')=\theta_j\,p_j(\theta_j x,\theta_j^{-1}\xi')\,,
$$
which is $x'$-frequency localized at scale $\mu^{\frac 12}$, 
where $\mu=\theta\l$ is the
frequency scale of the rescaled solution $u(\theta x)$ 
(we suppress the index $j$.)
We work with the wave packet transform $\tilde u$ of $u$ 
with respect to the $x'$ variables,
and let $\Theta$ denote the Hamiltonian flow along $\xi_1-q(x,\xi')$.
The reduction steps of sections 2 through 4 can then be adapted 
to reduce matters to
establishing the following.

\begin{theorem}\label{lastmaintheorem}
Suppose that $f\in L^2(\R^{2n})$ is supported in
a set of the form 
$\xi_{n+1}\approx\mu\,,$ $|\xi_2,\ldots,\xi_{n-1}|\le c\mu\,,$
$\xi_n\approx\theta\mu\,$
or
$|\xi_n|\lesssim\mu^{\frac 12}$ in case
$\theta=\mu^{-\frac 12}$.

Then, if
$u=T^*_\mu\bigl[f\bigl(\Theta_{0,x_1}(x',\xi')\bigr)\bigr]\,,$ we
have for $q\ge \frac{2n}{n-2}$
$$
\|u\|_{L^qL^2(S)}
\lesssim\mu^{\delta(q)}\theta^{\frac 12-\frac 1q}
\,\|f\|_{L^2(\R^{2n})}\,,
$$
and for $\frac{2(n+1)}{n-1}\le q\le \frac{2n}{n-2}$
$$
\|u\|_{L^qL^2(S)}
\lesssim\mu^{\delta(q)}\theta^{(n-1)(\frac 12-\frac 1q)-\frac 2q}
\,\|f\|_{L^2(\R^{2n})}\,.
$$
\end{theorem}

This implies Theorem \ref{penultimatemaintheorem} for $q$ 
such that the exponent of $\theta$
is at least $\frac 1q$. For $n\ge 4$, this happens for 
$q\ge 4\ge \frac {2n}{n-2}$. For $n=3$
this holds for $q\ge 5$.

In applying the reductions of section \ref{section2}, 
care must be taken since $\delta(q)\ge 1$
for large $n$, whereas the commutator $[A,\Gamma(D)]$ 
maps $H^{\delta-1}\rightarrow H^\delta$
only for $0\le\delta\le 1$. Here, $\Gamma(D)$ is a conic cutoff to the set
$|\xi_{n+1}|\approx |\xi_1,\ldots,\xi_n|\,.$ 
To get around this problem, in case $\delta(q)\ge 1$
we write $\delta(q)=m+\delta$, with $0\le\delta\le 1$. 
Let $d_T=(d_1,\ldots,d_{n-1},d_{n+1})$ denote the
tangential derivatives, and $d_T^m$ the collection of 
tangential derivatives of order at most $m$.
Then the extended and $\phi$-localized solution $u$ satisfies
$$
\|d_T^m u\|_{H^\delta}+\|d_T^m F\|_{H^\delta}\lesssim 
\|f\|_{H^{\delta(q)}(M)}\,.
$$
Since $d_T^m A$ is Lipschitz, it is easy to see that
$$
\|d_T^m[A,\Gamma(D)]Du\|_{H^\delta}\lesssim \|d_T^m u\|_{H^\delta}\,.
$$
We also gain powers of $d_T^m$ in the elliptic regularity arguments, 
and deduce that
$$
\|d_T^m(1-\Gamma(D))u\|_{H^{\delta+1}}\lesssim 
\|d_T^m u\|_{H^\delta}+\|d_T^m F\|_{H^\delta}\,.
$$
The norm on the left is sufficient to control 
$\|(1-\Gamma(D))u\|_{L^q_xL^2_t}$, and
we are reduced to considering $\Gamma(D)u$. 
This term, however, has Fourier transform
supported
%away from
outside of a conic neighborhood of
 the $\xi_n$ axis, hence
$$
\|\Gamma(D)u\|_{H^{\delta(q)}}\approx \|d_T^m \Gamma(D)u\|_{H^\delta}\,.
$$
The remaining reductions of section \ref{section2} then follow.

To prove Theorem \ref{lastmaintheorem}, we establish 
mapping properties for the kernel $K$ of $WW^*$,
localized in $\zeta=(\zeta_2,\ldots,\zeta_{n+1})$ 
by a cutoff $\beta_\theta(\zeta)$ to the set
$$
\zeta_{n+1}\approx \mu\,,\qquad |(\zeta_2,\ldots,\zeta_{n-1})|
\le c\,\mu\,,\qquad \zeta_n\approx \theta\mu\,,
$$
(respectively $|\zeta_n|\le \mu^{-\frac 12}$ in case $\theta=\mu^{-\frac 12}$.)
The bounds we establish, analogous to \eqref{kest1} and \eqref{kest2}, are
\begin{equation}\label{kest1'}
\sup_{r,s\in[0,\eps]}\Bigl\|\int K(r,x';s,y')\,f(y')\,dy'\Bigr\|_{L^2_{x'}}\le
\|f\|_{L^2_{y'}}\,,
\end{equation}
and
\begin{multline}\label{kest2'}
\Bigl\|\int K(r,x';s,y')\,f(y')\,dy
\Bigr\|_{L^\infty_{x_2,\ldots,x_n}L^2_{x_{n+1}}}
\\
\lesssim
\mu^{n-1}\theta\,\bigl(\,1+\mu\,|r-s|\,\bigr)^{-\frac{n-2}2}\,
\bigl(\,1+\mu\theta^2\,|r-s|\,\bigr)^{-\frac 12}\,
\|f\|_{L^1_{y_2,\ldots,y_n}L^2_{y_{n+1}}}\,.
\end{multline}
To see that this implies Theorem \ref{lastmaintheorem}, 
note that interpolation yields the bound
\begin{multline*}
\Bigl\|\int K(r,x';s,y')\,f(y')\,dy
\Bigr\|_{L^q_{x_2,\ldots,x_n}L^2_{x_{n+1}}}\\
\lesssim
\bigl(\mu^{n-1}\theta\bigr)^{1-\frac 2q}\,
\bigl(\,1+\mu|r-s|\,)^{-(n-2)(\frac 12-\frac 1q)}\,
\bigl(\,1+\mu\theta^2|r-s|\,\bigr)^{-(\frac 12-\frac 1q)}\,
\|f\|_{L^{q'}_{y_2,\ldots,y_n}L^2_{y_{n+1}}}\,.
\end{multline*}
If $q\ge \frac{2n}{n-2}\,,$ then $(n-2)(\frac 12-\frac 1q)\ge \frac 2q\,,$ 
and by the Hardy-Littlewood-Sobolev lemma
we obtain
$$
\Bigl\|\int K(r,x';s,y')\,F(s,y')\,dy
\Bigr\|_{L^q_{r,x_2,\ldots,x_n}L^2_{x_{n+1}}}
\lesssim
\mu^{2\delta(q)}\theta^{1-\frac 2q}\,
\|F\|_{L^{q'}_{s,y_2,\ldots,y_n}L^2_{y_{n+1}}}\,.
$$
If $\frac{2n}{n-2}\ge q\ge\frac{2(n+1)}{n-1}\,,$ then
$$
\bigl(\,1+\mu|r-s|\,)^{-(n-2)(\frac 12-\frac 1q)}\,
\bigl(\,1+\mu\theta^2|r-s|\,\bigr)^{-(\frac 12-\frac 1q)}\le
\mu^{-\frac 2q}\theta^{-\frac 4q+(n-2)(1-\frac 2q)}\,|r-s|^{-\frac 2q}\,,
$$
and Hardy-Littlewood-Sobolev yields Theorem \ref{lastmaintheorem} 
for this case.

We now turn to the proof of \eqref{kest1'} and \eqref{kest2'}. 
The estimate \eqref{kest1'}
follows as does the estimate \eqref{kest1} from the 
boundedness of $T_\mu$ and the fact that $\Theta_{r,s}$
preserves the measure $dx'\,d\xi'$.
To establish \eqref{kest2'}, we consider as before separate cases, 
depending on $|r-s|$.

Consider the case $\mu\theta^2|r-s|\ge 1$. We fix $\thetabar\le\theta$ 
so that $\mu\thetabar^2|r-s|=1$,
and decompose $\beta_\theta(\zeta)$ into a sum of cutoffs $\beta_j(\zeta)$, 
each of which
is localized to a cone of angle $\thetabar$ about some direction. 
As in the proof of \eqref{kest3},
we have that
$$
\int|K_j(r,x';s,y')|\,dy_{n+1}\,\lesssim\,
\mu^{n-1}\thetabar^{n-1}
\bigl(\,1+\mu\thetabar|(y'-w^j_{s,r})_{2,\ldots,n}|\bigr)^{-N}\,,
$$
where the $w^j_{s,r}$ give a $(\mu\thetabar)^{-1}$ 
separated set after projection onto the $(2,\ldots,n)$ variables.
Adding over $j$ yields the desired bounds, since
$$
\mu^{n-1}\thetabar^{n-1}=\mu^{\frac {n-1}2}\,|r-s|^{-\frac{n-1}2}\,.
$$

In case $\mu\theta^2|r-s|\le 1$, let $\thetabar\ge\theta$ be given by
$$
\thetabar=\min\bigl(\,\mu^{-\frac 12}|r-s|^{-\frac 12}\,,\,1\,\bigr)\,.
$$
We set $\zeta''=(\zeta_2,\ldots,\zeta_{n-1},\zeta_{n+1})$,
and let $\beta_j$ be a partition of unity in cones of angle 
$\thetabar$ on $\R^{n-1}$.
We then decompose
$$
\beta_\theta(\zeta)=\sum_j \beta_\theta(\zeta)\,\beta_j(\zeta'')\,.
$$
Let $K=\sum_j K_j$ denote the corresponding kernel decomposition. 
As in the proof of Theorem \ref{ktheorem},
we can bound $K_j$ by
\begin{multline*}
\mu^{\frac n2}\int\bigl(\,1+\mu\thetabar\,|d_{\zeta''} \zeta_{s,r}
\cdot(y'-x'_{s,r})|+
\mu\theta\,|d_{\zeta_n} \zeta_{s,r}\cdot(y'-x'_{s,r})|+
|\langle\zeta_{s,r},y'-x'_{s,r}\rangle|
\,\bigr)^{-N}\\
\times\bigl(\,1+\mu^{\frac 12}|x'-z|\,\bigr)^{-N}\,dz\,d\zeta\,.
\end{multline*}
Here, $(x'_{s,r},\xi'_{s,r})=\Theta_{s,r}(x',\xi'_j)$, with $\xi'_j$
a fixed vector in the support of $\beta_\theta(\zeta)\beta_j(\zeta'')$. 
Also, $(z_{s,r},\zeta_{s,r})=\Theta_{s,r}(z,\zeta)$.
Since $d_\zeta \zeta_{s,r}$ is invertible, and 
$\mu\thetabar\ge\mu\theta\ge\mu^{\frac 12}$,
the first two terms in the integrand dominate $\mu^{\frac 12}|y'-x'_{s,r}|$.

We first show that we may replace $\zeta_{s,r}$ by 
$\xi'_{s,r}=\zeta_{s,r}(x',\xi'_j)$ in the third term in
parentheses above. By homogeneity, we may consider $|\zeta|=|\xi'_j|$.
We take a first order Taylor expansion, and use bounds 
\eqref{esttwo} on $d_\zeta^2\zeta_{s,r}$,
to write
$$
\zeta_{s,r}-\zeta_{s,r}(z,\xi'_j)=(\zeta-\xi'_j)\cdot d_\zeta\zeta_{s,r}+
O(\,|\zeta-\xi'_j|^2\mu^{-\frac 12}|s-r|\,)\,.
$$
Since
$$
|(\zeta-\xi'_j)''|\lesssim \mu\thetabar\,,\qquad
|(\zeta-\xi'_j)_n|\lesssim \mu\theta\,,\qquad
\mu\thetabar^2|s-r|\le 1\,,
$$
this shows we may replace $\zeta_{s,r}$ by $\zeta_{s,r}(z,\xi'_j)$, 
as the errors are
absorbed by the first two terms in parentheses.
On the other hand, by \eqref{estone}
$$
|\,\langle\zeta_{s,r}(x',\xi'_j)-\zeta(z,\xi'_j)\,,y'-x'_{s,r}\rangle\,|
\lesssim 
\mu\,|x'-z|\,|y'-x'_{s,r}|\,,
$$
which is also absorbed by the other terms.

We next use \eqref{estone} to see that we may replace 
$d_\zeta\zeta_{s,r}$ by the identity
matrix, since the error induced is dominated by
$$
\mu\thetabar|s-r|\,|y'-x'_{s,r}|\le \mu^{\frac 12}|y'-x'_{s,r}|\,.
$$

Consequently, since $\xi'_{s,r}$ has $n+1$ component comparable to 
$\mu$, we obtain
$$
\int|K_j(r,x';s,y')|\,dy_{n+1}
\lesssim\mu^{n-1}\thetabar^{n-2}\,\theta\,
\bigl(\,1+\mu\thetabar\,|(y'-x'_{s,r})_{2,\ldots,n-1}|\,\bigr)^{-N}\,.
$$
The points $x'_{s,r}$ are $\mu\thetabar$ separated in the 
$(2,\ldots,n-1)$ variables as $j$ varies,
which follows by Corollary \ref{dzetazest} and the fact that 
$q(z,\zeta)$ is close to $|\zeta|$,
hence we can add over $j$ to obtain
$$
\int|K(r,x';s,y')|\,dy_{n+1}
\lesssim \mu^{n-1}\thetabar^{n-2}\,\theta\lesssim
\mu^{n-1}\theta\bigl(\,1+\mu|r-s|\,\bigr)^{-\frac{n-2}2}\,.
\qed
$$

\end{document}